\documentclass[12pt,reqno]{amsart}

\usepackage{amscd}
\usepackage{epsfig}
\usepackage{float}
\usepackage{url}
\usepackage[all]{xy}

\newcommand{\N}{{\mathbb N}}
\newcommand{\Z}{{\mathbb Z}}
\newcommand{\Q}{{\mathbb Q}}
\newcommand{\C}{{\mathbb C}}
\newcommand{\R}{{\mathbb R}}

\newcommand{\OO}{{\mathcal O}}

\newcommand{\www}{\widetilde}
\newcommand{\wwh}{\widehat}
\newcommand{\oooo}{\overline}

\newcommand{\paa}{\partial}

\DeclareMathOperator{\divis}{div}

\DeclareMathOperator{\id}{id}

\DeclareMathOperator{\lcm}{lcm}

\DeclareMathOperator{\supp}{supp}

\newlength{\unitlengthGrB}
\newlength{\unitlengthGrC}
\newlength{\unitlengthGrD}
\newcommand{\GrBi}[1][1]
{
	\xygraph{
		!{<0cm,0cm>;<#1\unitlengthGrB,0cm>:<0cm,#1\unitlengthGrB>::}
		!{(0,0) }*+{\bullet_{1}}="1"
		!{(2,0) }*+{\bullet_{2}}="2"
	}
}

\newcommand{\GrBii}[1][1]
{
	\xygraph{
		!{<0cm,0cm>;<#1\unitlengthGrB,0cm>:<0cm,#1\unitlengthGrB>::}
		!{(0,0) }*+{\bullet_{1}}="1"
		!{(2,0) }*+{\bullet_{2}}="2"
		"1":"2"
	}
}

\newcommand{\GrBiii}[1][1]
{
	\xygraph{
		!{<0cm,0cm>;<#1\unitlengthGrB,0cm>:<0cm,#1\unitlengthGrB>::}
		!{(0,0) }*+{\bullet_{1}}="1"
		!{(2,0) }*+{\bullet_{2}}="2"
		"1":@/_/"2" "2":@/_/"1"
	}
}

\newcommand{\GrCi}[1][1]
{
	\xygraph{
		!{<0cm,0cm>;<#1\unitlengthGrC,0cm>:<0cm,#1\unitlengthGrC>::}
		!{(0,0) }*+{\bullet_{1}}="1"
		!{(2,0) }*+{\bullet_{2}}="2"
		!{(1,1.7) }*+{\bullet_{3}}="3"
	}
}

\newcommand{\GrCii}[1][1]
{
	\xygraph{
		!{<0cm,0cm>;<#1\unitlengthGrC,0cm>:<0cm,#1\unitlengthGrC>::}
		!{(0,0) }*+{\bullet_{1}}="1"
		!{(2,0) }*+{\bullet_{2}}="2"
		!{(1,1.7) }*+{\bullet_{3}}="3"
		"3":"2"
	}
}

\newcommand{\GrCiii}[1][1]
{
	\xygraph{
		!{<0cm,0cm>;<#1\unitlengthGrC,0cm>:<0cm,#1\unitlengthGrC>::}
		!{(0,0) }*+{\bullet_{1}}="1"
		!{(2,0) }*+{\bullet_{2}}="2"
		!{(1,1.7) }*+{\bullet_{3}}="3"
		"2":"1" "3":"1"
	}
}

\newcommand{\GrCiv}[1][1]
{
	\xygraph{
		!{<0cm,0cm>;<#1\unitlengthGrC,0cm>:<0cm,#1\unitlengthGrC>::}
		!{(0,0) }*+{\bullet_{1}}="1"
		!{(2,0) }*+{\bullet_{2}}="2"
		!{(1,1.7) }*+{\bullet_{3}}="3"
		"2":@/_/"3" "3":@/_/"2"
	}
}

\newcommand{\GrCvi}[1][1]
{
	\xygraph{
		!{<0cm,0cm>;<#1\unitlengthGrC,0cm>:<0cm,#1\unitlengthGrC>::}
		!{(0,0) }*+{\bullet_{1}}="1"
		!{(2,0) }*+{\bullet_{2}}="2"
		!{(1,1.7) }*+{\bullet_{3}}="3"
		"1":@/_/"2" "2":@/_/"1" "3":"1"
	}
}

\newcommand{\GrCv}[1][1]
{
	\xygraph{
		!{<0cm,0cm>;<#1\unitlengthGrC,0cm>:<0cm,#1\unitlengthGrC>::}
		!{(0,0) }*+{\bullet_{1}}="1"
		!{(2,0) }*+{\bullet_{2}}="2"
		!{(1,1.7) }*+{\bullet_{3}}="3"
		"2":"1" "3":"2"
	}
}

\newcommand{\GrCvii}[1][1]
{
	\xygraph{
		!{<0cm,0cm>;<#1\unitlengthGrC,0cm>:<0cm,#1\unitlengthGrC>::}
		!{(0,0) }*+{\bullet_{1}}="1"
		!{(2,0) }*+{\bullet_{2}}="2"
		!{(1,1.7) }*+{\bullet_{3}}="3"
		"1":"2" "2":"3" "3":"1"
	}
}

\newcommand{\GrDi}[1][1]
{
	\xygraph{
		!{<0cm,0cm>;<#1\unitlengthGrD,0cm>:<0cm,#1\unitlengthGrD>::}
		!{(0,0) }*+{\bullet_{1}}="1"
		!{(2,0) }*+{\bullet_{2}}="2"
		!{(2,2) }*+{\bullet_{3}}="3"
		!{(0,2) }*+{\bullet_{4}}="4"
	}
}

\newcommand{\GrDiii}[1][1]
{
	\xygraph{
		!{<0cm,0cm>;<#1\unitlengthGrD,0cm>:<0cm,#1\unitlengthGrD>::}
		!{(0,0) }*+{\bullet_{1}}="1"
		!{(2,0) }*+{\bullet_{2}}="2"
		!{(2,2) }*+{\bullet_{3}}="3"
		!{(0,2) }*+{\bullet_{4}}="4"
		"3":@/_/"4" "4":@/_/"3"
	}
}

\newcommand{\GrDxviii}[1][1]
{
	\xygraph{
		!{<0cm,0cm>;<#1\unitlengthGrD,0cm>:<0cm,#1\unitlengthGrD>::}
		!{(0,0) }*+{\bullet_{1}}="1"
		!{(2,0) }*+{\bullet_{2}}="2"
		!{(2,2) }*+{\bullet_{3}}="3"
		!{(0,2) }*+{\bullet_{4}}="4"
		"1":@/_/"2" "2":@/_/"1" "3":@/_/"4" "4":@/_/"3"
	}
}

\newcommand{\GrDvii}[1][1]
{
	\xygraph{
		!{<0cm,0cm>;<#1\unitlengthGrD,0cm>:<0cm,#1\unitlengthGrD>::}
		!{(0,0) }*+{\bullet_{1}}="1"
		!{(2,0) }*+{\bullet_{2}}="2"
		!{(2,2) }*+{\bullet_{3}}="3"
		!{(0,2) }*+{\bullet_{4}}="4"
		"2":"3" "3":"4" "4":"2"
	}
}

\newcommand{\GrDxiv}[1][1]
{
	\xygraph{
		!{<0cm,0cm>;<#1\unitlengthGrD,0cm>:<0cm,#1\unitlengthGrD>::}
		!{(0,0) }*+{\bullet_{1}}="1"
		!{(2,0) }*+{\bullet_{2}}="2"
		!{(2,2) }*+{\bullet_{3}}="3"
		!{(0,2) }*+{\bullet_{4}}="4"
		"1":"2" "2":"3" "3":"4" "4":"1"
	}
}

\newcommand{\GrDii}[1][1]
{
	\xygraph{
		!{<0cm,0cm>;<#1\unitlengthGrD,0cm>:<0cm,#1\unitlengthGrD>::}
		!{(0,0) }*+{\bullet_{1}}="1"
		!{(2,0) }*+{\bullet_{2}}="2"
		!{(2,2) }*+{\bullet_{3}}="3"
		!{(0,2) }*+{\bullet_{4}}="4"
		"4":"3"
	}
}

\newcommand{\GrDviii}[1][1]
{
	\xygraph{
		!{<0cm,0cm>;<#1\unitlengthGrD,0cm>:<0cm,#1\unitlengthGrD>::}
		!{(0,0) }*+{\bullet_{1}}="1"
		!{(2,0) }*+{\bullet_{2}}="2"
		!{(2,2) }*+{\bullet_{3}}="3"
		!{(0,2) }*+{\bullet_{4}}="4"
		"2":@/_/"3" "3":@/_/"2" "4":"3"
	}
}

\newcommand{\GrDix}[1][1]
{
	\xygraph{
		!{<0cm,0cm>;<#1\unitlengthGrD,0cm>:<0cm,#1\unitlengthGrD>::}
		!{(0,0) }*+{\bullet_{1}}="1"
		!{(2,0) }*+{\bullet_{2}}="2"
		!{(2,2) }*+{\bullet_{3}}="3"
		!{(0,2) }*+{\bullet_{4}}="4"
		"1":@/_/"2" "2":@/_/"1" "4":"3"
	}
}

\newcommand{\GrDxv}[1][1]
{
	\xygraph{
		!{<0cm,0cm>;<#1\unitlengthGrD,0cm>:<0cm,#1\unitlengthGrD>::}
		!{(0,0) }*+{\bullet_{1}}="1"
		!{(2,0) }*+{\bullet_{2}}="2"
		!{(2,2) }*+{\bullet_{3}}="3"
		!{(0,2) }*+{\bullet_{4}}="4"
		"1":"2" "2":"3" "3":"1" "4":"3"
	}
}

\newcommand{\GrDvi}[1][1]
{
	\xygraph{
		!{<0cm,0cm>;<#1\unitlengthGrD,0cm>:<0cm,#1\unitlengthGrD>::}
		!{(0,0) }*+{\bullet_{1}}="1"
		!{(2,0) }*+{\bullet_{2}}="2"
		!{(2,2) }*+{\bullet_{3}}="3"
		!{(0,2) }*+{\bullet_{4}}="4"
		"3":"2" "4":"3"
	}
}

\newcommand{\GrDxvi}[1][1]
{
	\xygraph{
		!{<0cm,0cm>;<#1\unitlengthGrD,0cm>:<0cm,#1\unitlengthGrD>::}
		!{(0,0) }*+{\bullet_{1}}="1"
		!{(2,0) }*+{\bullet_{2}}="2"
		!{(2,2) }*+{\bullet_{3}}="3"
		!{(0,2) }*+{\bullet_{4}}="4"
		"1":@/_/"2" "2":@/_/"1" "3":"2" "4":"3"
	}
}

\newcommand{\GrDx}[1][1]
{
	\xygraph{
		!{<0cm,0cm>;<#1\unitlengthGrD,0cm>:<0cm,#1\unitlengthGrD>::}
		!{(0,0) }*+{\bullet_{1}}="1"
		!{(2,0) }*+{\bullet_{2}}="2"
		!{(2,2) }*+{\bullet_{3}}="3"
		!{(0,2) }*+{\bullet_{4}}="4"
		"2":"1" "3":"2" "4":"3"
	}
}

\newcommand{\GrDiv}[1][1]
{
	\xygraph{
		!{<0cm,0cm>;<#1\unitlengthGrD,0cm>:<0cm,#1\unitlengthGrD>::}
		!{(0,0) }*+{\bullet_{1}}="1"
		!{(2,0) }*+{\bullet_{2}}="2"
		!{(2,2) }*+{\bullet_{3}}="3"
		!{(0,2) }*+{\bullet_{4}}="4"
		"2":"1" "4":"3"
	}
}

\newcommand{\GrDv}[1][1]
{
	\xygraph{
		!{<0cm,0cm>;<#1\unitlengthGrD,0cm>:<0cm,#1\unitlengthGrD>::}
		!{(0,0) }*+{\bullet_{1}}="1"
		!{(2,0) }*+{\bullet_{2}}="2"
		!{(2,2) }*+{\bullet_{3}}="3"
		!{(0,2) }*+{\bullet_{4}}="4"
		"3":"2" "4":"2"
	}
}

\newcommand{\GrDxvii}[1][1]
{
	\xygraph{
		!{<0cm,0cm>;<#1\unitlengthGrD,0cm>:<0cm,#1\unitlengthGrD>::}
		!{(0,0) }*+{\bullet_{1}}="1"
		!{(2,0) }*+{\bullet_{2}}="2"
		!{(2,2) }*+{\bullet_{3}}="3"
		!{(0,2) }*+{\bullet_{4}}="4"
		"1":@/_/"2" "2":@/_/"1" "3":"2" "4":"1"
	}
}

\newcommand{\GrDxix}[1][1]
{
	\xygraph{
		!{<0cm,0cm>;<#1\unitlengthGrD,0cm>:<0cm,#1\unitlengthGrD>::}
		!{(0,0) }*+{\bullet_{1}}="1"
		!{(2,0) }*+{\bullet_{2}}="2"
		!{(2,2) }*+{\bullet_{3}}="3"
		!{(0,2) }*+{\bullet_{4}}="4"
		"1":@/_/"2" "2":@/_/"1" "3":"1" "4":"1"
	}
}

\newcommand{\GrDxii}[1][1]
{
	\xygraph{
		!{<0cm,0cm>;<#1\unitlengthGrD,0cm>:<0cm,#1\unitlengthGrD>::}
		!{(0,0) }*+{\bullet_{1}}="1"
		!{(2,0) }*+{\bullet_{2}}="2"
		!{(2,2) }*+{\bullet_{3}}="3"
		!{(0,2) }*+{\bullet_{4}}="4"
		"2":"1" "3":"2" "4":"1"
	}
}

\newcommand{\GrDxi}[1][1]
{
	\xygraph{
		!{<0cm,0cm>;<#1\unitlengthGrD,0cm>:<0cm,#1\unitlengthGrD>::}
		!{(0,0) }*+{\bullet_{1}}="1"
		!{(2,0) }*+{\bullet_{2}}="2"
		!{(2,2) }*+{\bullet_{3}}="3"
		!{(0,2) }*+{\bullet_{4}}="4"
		"2":"1" "3":"2" "4":"2"
	}
}

\newcommand{\GrDxiii}[1][1]
{
	\xygraph{
		!{<0cm,0cm>;<#1\unitlengthGrD,0cm>:<0cm,#1\unitlengthGrD>::}
		!{(0,0) }*+{\bullet_{1}}="1"
		!{(2,0) }*+{\bullet_{2}}="2"
		!{(2,2) }*+{\bullet_{3}}="3"
		!{(0,2) }*+{\bullet_{4}}="4"
		"2":"1" "3":"1" "4":"1"
	}
}

\begin{document}

\theoremstyle{plain}
\newtheorem{lemma}{Lemma}[section]
\newtheorem{definition/lemma}[lemma]{Definition/Lemma}
\newtheorem{theorem}[lemma]{Theorem}
\newtheorem{proposition}[lemma]{Proposition}
\newtheorem{corollary}[lemma]{Corollary}
\newtheorem{conjecture}[lemma]{Conjecture}
\newtheorem{conjectures}[lemma]{Conjectures}

\theoremstyle{definition}
\newtheorem{definition}[lemma]{Definition}
\newtheorem{withouttitle}[lemma]{}
\newtheorem{remark}[lemma]{Remark}
\newtheorem{remarks}[lemma]{Remarks}
\newtheorem{example}[lemma]{Example}
\newtheorem{examples}[lemma]{Examples}

\title[Classification of quasihomogeneous singularities]
{On the classification of quasihomogeneous singularities} 

\author{Claus Hertling \and Ralf Kurbel}

\address{Claus Hertling\\
Lehrstuhl f\"ur Mathematik VI, Universit\"at Mannheim, Seminargeb\"aude
A 5, 6, 68131 Mannheim, Germany}

\email{hertling@math.uni-mannheim.de}

\address{Ralf Kurbel\\
Lehrstuhl f\"ur Mathematik VI, Universit\"at Mannheim, Seminargeb\"aude
A 5, 6, 68131 Mannheim, Germany}

\email{kurbel@math.uni-mannheim.de}





\begin{abstract}
{\small The motivation for this paper are computer calculations of complete
lists of weight systems of quasihomogeneous polynomials with 
isolated singularity at $0$ up to rather large Milnor numbers.
We review combinatorial characterizations of such weight systems
for any number of variables.
This leads to certain types and graphs of such weight systems.
Using them,
we prove an upper bound for the common denominator 
(and the order of the monodromy) by the Milnor number, 
and we show surprising consequences if the Milnor number is a prime number.
}
\end{abstract}

\maketitle

\tableofcontents

\setcounter{section}{0}

\section{Introduction}\label{c1}
\setcounter{equation}{0}

\noindent
Several people have (re)discovered characterizations of those weight systems
which admit quasihomogeneous polynomials with isolated singularity at $0$.
Section \ref{s2} collects and compares these characterizations and gives
all references which we found. The results of this section are not new.
But the references are not well known and for several reasons, it is not
so easy to extract the results from them.
Also, we will need part of the characterizations 
for a good control of such weight systems in the later sections.

In section \ref{s3} a part of the conditions
is used to associate after a choice a type and a graph to a quasihomogeneous
singularity. The idea for this is contained in \cite{Ar}\cite{AGV}, and there it is carried out
for $2$ and $3$ variables. The general case is carried out in \cite{OR1},
but that part of \cite{OR1} was never published. As we will need the graphs
in the sections \ref{s4} and \ref{s6}, we rewrite the general case.
Section \ref{s3} also makes the classification in the case of $4$ variables
in \cite{YS} more precise, showing how necessary and sufficient conditions are obtained.

Section \ref{s4} gives an estimate $d\leq \textup{const}(n)\cdot \mu$
for the weighted degree $d$ of a reduced weight system of a quasihomogeneous singularity
from above by the Milnor number $\mu$.
The calculations start with the well known formula for $\mu$ in terms of the
weights, but refine this formula using a graph and a type of the singularity.
The estimate is useful for a computer calculation of all reduced weight systems
of quasihomogeneous singularities up to a given Milnor number. 
We carried out such computer calculations for $2$, $3$ and $4$ variables and 
$\mu\leq 9000$, $9000$ and $2000$.
The long tables are available on the homepage \cite{HK}. Some observations from them are
formulated in section \ref{s5}.

Section \ref{s6} proves a surprising fact which we found looking at these tables.
If the Milnor number of a quasihomogeneous singularity is a prime number,
then the only type which one can associate to it is the chain type
(up to adding or removing squares from the singularity), 
and furthermore, all eigenvalues of the monodromy have multiplicity one. 
The proof further refines the formula for the Milnor number from section \ref{s4}.

We thank Sasha Aleksandrov for translating \cite{Ko2} and for the reference 
\cite{YS} and Wolfgang Ebeling and 
Atsushi Takahashi for discussions related to lemma \ref{t3.5} and \cite{ET}.
This paper was written during a stay at the Tokyo Metropolitan University.
We thank the TMU and Martin Guest for hospitality.

\section{Combinatorial characterizations of weight systems 
of quasihomogeneous singularities}\label{s2}
\setcounter{equation}{0}

\noindent
We note $\N_0=\{0,1,2,...\}\supset \N=\{1,2,...\}$. The support of a polynomial

$$f=\sum_{\alpha\in\N_0^n}a_\alpha\cdot x^\alpha\in \C[x_1,...,x_n]
\quad \textup{where}\quad x^\alpha=x_1^{\alpha_1}...x_n^{\alpha_n}$$ 
is $\supp f=\{\alpha\in\N_0^n\ |\ a_\alpha\neq 0\}$.
The polynomial is called {\it quasihomogeneous} with weight system
$(w_1,...,w_n,d)\in\R_{>0}^{n+1}$ if
$$\sum_{i=1}^n\alpha_i\cdot w_i=d\quad \textup{ for any }\alpha\in\supp f.$$
Here $w_1,...,w_n$ are the weights and $d$ is the weighted degree.
If a polynomial is quasihomogeneous with some weight system it is also quasihomogeneous
with a weight system $(w_1,...,w_n,d)\in\Q_{>0}^{n+1}$.
If a quasihomogeneous polynomial has an isolated singularity at $0$, that is,
if the $\frac{\paa f}{\paa x_i}$ vanish simultaneously precisely at $0$, then
$w_i<d$ for all $i$. 
Therefore, from now on throughout the whole paper we consider only weight systems
with
$$(w_1,...,w_n,d)\in\Q_{>0}^{n+1}\qquad \textup{ and }\qquad 
w_i<d\textup{ for all }i.$$
Furthermore, from now on we reserve the letters $v_1,...,v_n$ for weights of weight systems
$(v_1,...,v_n,d)\in\N^{n+1}$, and the letters $w_1,...,w_n$ for weights of normalized 
weight systems $(w_1,...,w_n,1)\in\Q_{>0}^{n+1}$, that is, with weighted degree 1.

A weight system $(v_1,...,v_n,d)\in\N^{n+1}$ is called {\it reduced} 
if $\gcd(v_1,...,v_n,d)=1$. 
In later chapters, but not in this one, we will also use a result in \cite{Sa1}
and restrict to weight systems with $v_i\leq \frac{d}{2}$ and $w_i\leq\frac{1}{2}$.

Fix $n\in\N$ and denote $N:=\{1,...,n\}$ and 
$e_i:=(\delta_{ij})_{j=1,...,n}\in\N_0^n$.
For $J\subset N$ and a weight system $(v_1,...,v_n,d)\in\N^{n+1}$ (with $v_i<d$) and 
$k\in\N_0$ denote
\begin{eqnarray*}
\N_0^J&:=& \{\alpha\in\N_0^n\ |\ \alpha_i=0\textup{ for }i\notin J\},\\
(\N_0^n)_k&:=& \{\alpha\in \N_0^n\ |\ \sum_i\alpha_i\cdot v_i=k\},\\
(\N_0^J)_k&:=& \N_0^J\cap (\N_0^n)_k.
\end{eqnarray*}
The following combinatorial lemma will help to compare in theorem \ref{t2.2}
several characterizations of weight systems which admit quasihomogeneous polynomials
with isolated singularities. A discussion of the history and references will
be given after theorem \ref{t2.2}.

\begin{lemma}\label{t2.1}
Fix a weight system $(v_1,...,v_n,d)\in\N^{n+1}$ with $v_i<d$ and a subset 
$R\subset (\N_0^n)_d$. For any $k\in N$ define
$$R_k:= \{\alpha\in (\N_0^n)_{d-v_k}\ |\ \alpha+e_k\in R\}.$$
The following five conditions (C1), (C1)', (C2), (C2)' and (C3) are equivalent.
\begin{list}{}{}
\item[(C1):] 
\quad $\forall\ J\subset N\textup{ with }J\neq \emptyset$\\
\hspace*{2cm} $\exists\ \alpha\in R\cap \N_0^J$\\
\hspace*{1.6cm}or $\exists\ K\subset N-J\textup{ with }|K|=|J|$\\
\hspace*{5cm}$\textup{ and }\forall\ k\in K\ \exists\ \alpha\in R_k\cap \N_0^J.$\\
\item[(C1)':] 
\quad As (C1), but only $J$ with $|J|\leq \frac{n+1}{2}$.\\
\item[(C2):] 
\quad $\forall\ J\subset N\textup{ with }J\neq \emptyset$\\
\hspace*{2cm}$\exists\ K\subset N\textup{ with }|K|=|J|$\\
\hspace*{4cm}\textup{ and }$\forall\ k\in K\ \exists\ \alpha\in R_k\cap \N_0^J.$\\
\item[(C2)':] 
\quad As (C2), but only $J$ with $|J|\leq \frac{n+1}{2}$.\\
\item[(C3):] 
\quad $\forall\ I,\ J\subset N\textup{ with }|I|<|J|$\\
\hspace*{2cm}$\exists\ k\in N-I\textup{ and }\exists\ \alpha\in R_k\cap \N_0^J.$
\end{list}
\end{lemma}

{\bf Proof:}
(C1)$\Rightarrow$(C1)' and (C2)$\Rightarrow$(C2)' are trivial.

(C1)'$\Rightarrow$(C1): Consider $J\subset N$ with $|J|> \frac{n+1}{2}$
and $I\subset J$ with $|I|=\left[\frac{n+1}{2}\right]$.
If there exists $\alpha\in R\cap \N_0^I$, then also $\alpha\in R\cap\N_0^J$.
If not, then there exists $K\subset N-I$ with $|K|=|I|$ and
$\forall\ k\in K\ \exists\ \alpha\in R_k\cap \N_0^I$.
Then $n$ is even and $K=N-I$, and $K\cap J\neq\emptyset$, and
for $k\in K\cap J$ and $\alpha\in R_k\cap \N_0^I$ one finds
$\alpha+e_k\in R\cap \N_0^J$.

(C2)'$\Rightarrow$(C1)': Consider $J\subset N$ with $0<|J|\leq \frac{n+1}{2}$
and $K\subset N$ such that $J$ and $K$ satisfy (C2)'.

{\bf 1st case:} $K\subset N-J$. Then $J$ and $K$ satisfy (C1)'.

{\bf 2nd case:} $K\cap J\neq\emptyset$. Then for $k\in K\cap J$ and 
$\alpha\in R_k\cap \N_0^J$ one obtains $\alpha+e_k\in R\cap \N_0^J$,
so $J$ satisfies (C1)'.

(C3)$\Rightarrow$(C2): Consider $J\subset N$ with $J\neq \emptyset$. 
Construct elements $k_1,...,k_{|J|}\in N$ and subsets $I_j=\{k_1,...,k_j\}$
for $0\leq j\leq |J|-1$ and $K:=\{k_1,...,k_{|J|}\}$ as follows.
(C3) gives for $J$ and $I_j$ an element $k_{j+1}\in N-I_j$ with 
$R_{k_{j+1}}\cap \N_0^J\neq \emptyset$. Obviously $|K|=|J|$,
and $J$ and $K$ satisfy (C2).

(C1)$\Rightarrow$(C3): Consider $I,J\subset N$ with $|I|<|J|$. Then $J\neq \emptyset$.

{\bf 1st case}, $J$ and some set $K$ satisfy (C1): Because of $|I|<|J|=|K|$ there is a 
$k\in (N-I)\cap K$ with $R_k\cap \N_0^J\neq \emptyset$.

{\bf 2nd case}, $R\cap \N_0^J\neq \emptyset$: If $R\cap \N_0^J\neq R\cap \N_0^{J\cap I}$
then there exists an $\alpha\in R\cap \N_0^J-R\cap \N_0^{J\cap I}$ and
a $k\in J-J\cap I$ with $\alpha_k>0$. 
Then $k\in N-I$ and $\alpha-e_k\in R_k\cap \N_0^J$.

So suppose $R\cap \N_0^J=R\cap \N_0^{J\cap I}$. Then $J_1:= J-J\cap I\neq \emptyset$ 
because of $|I|<|J|$, and $R\cap \N_0^{J_1}=\emptyset$,
so there exists a $K_1\subset N-J_1$ such that $J_1$ and $K_1$ satisfy (C1).
If $K_1\cap J\neq \emptyset$ then for $k\in K_1\cap J$ and $\alpha\in R_k\cap \N_0^{J_1}$
one has $\alpha+e_k\in R\cap \N_0^J- R\cap \N_0^{J\cap I}$, a contradiction.
Thus $K_1\cap J=\emptyset$.

This and $|K_1|=|J_1|>|I-(J\cap I)|$ give $K_1-I=K_1-(I-(J\cap I))\neq\emptyset$.
Any $k\in K_1-I$ satisfies $\emptyset \neq R_k\cap \N_0^{J_1}\subset R_k\cap \N_0^J$.
\hfill$\Box$

\begin{theorem}\label{t2.2}
Let $(v_1,...,v_n,d)\in\N^{n+1}$ with $v_i<d$ be a weight system.

(a) Let $f\in\C[x_0,...,x_n]$ be a quasihomogeneous polynomial.
The condition
\begin{list}{}{}
\item[(IS1):] 
\quad $f$ has an isolated singularity at $0$,
\end{list}
implies that $R:=\supp f\subset (\N_0^n)_d$ satisfies (C1) to (C3).

\medskip
(b) Let $R$ be a subset of $(\N_0^n)_d$. The following conditions are equivalent.
\begin{list}{}{}
\item[(IS2):] \quad There exists a quasihomogeneous polynomial $f$ with $\supp f\subset R$
and an isolated singularity at 0.
\item[(IS2)':] \quad \quad A generic quasihomogeneous polynomial with $\supp f\subset R$ has an 
isolated singularity at 0.
\item[(C1) to (C3):] \quad $R$ satisfies (C1) to (C3).
\end{list}

\medskip
(c) In the case $R=(\N_0^n)_d$ obviously $R_k=(\N_0^n)_{d-v_k}$. 
The following conditions are equivalent.
\begin{list}{}{}
\item[(IS3):] \quad There exists a quasihomogeneous polynomial $f$ with
an isolated singularity at 0.
\item[(IS3)':] \quad A generic quasihomogeneous polynomial has an 
isolated singularity at 0.
\item[(C1) to (C3):] \quad $R=(\N_0^n)_d$ satisfies (C1) to (C3).
\end{list}
\end{theorem}

\begin{remarks}\label{t2.3}
Several people (re)discovered parts of this theorem. We will not reprove it here,
but comment on the history and the references.

(i) Of course, (IS2)$\iff$(IS2)' and (IS3)$\iff$(IS3)' and (b)$\Rightarrow$(c)
and (a)$\Rightarrow\bigl( \textup{(IS2)}\Rightarrow \textup{(C1) to (C3)}\bigr)$.

(ii) Part (a) is quite elementary, for example (IS1)$\Rightarrow$(C1) is contained
in K. Saito's paper \cite[Lemma 1.5]{Sa1}, and it can also be extracted
from \cite[Remark 3]{Sh}.

(iii) (IS2)$\iff$(C2) is part of an equivalence for more general functions in
\cite[Remarque 1.13 (ii)]{Ko1}, but there Kouchnirenko did not carry out the proof
in detail. He gave a short proof of the refined version (IS2)$\iff$(C2)'
in \cite[Theorem 1]{Ko2}. This reference \cite{Ko2} seems to have been cited 
up to now only in \cite{Sh}, it seems to have been almost completely ignored.

(iv) Around the same time as Kouchnirenko, Orlik and Randell proved (IS3)$\iff$(C3)
in the preprint \cite[Theorem 2.12]{OR1}, but the published part \cite{OR2} of it
does not contain this result. It seems that they have not published this result.

(v) O.P. Shcherbak stated a more general result \cite[Theorem 1]{Sh} from which
one can extract (IS2)$\iff$(C1), but he did not provide a proof.
That was done by Wall \cite[Ch. 5]{Wa}, who also stated explicitly
(IS2)$\iff$(C1) and (IS3)$\iff$(C1), they are Theorem 5-1 and Theorem 5-3
in \cite{Wa} for the hypersurface case (explicit in (5-7)).
But as he covers a much more general case, his proof is long.

(vi) A short proof of (IS3)$\iff$(C1) is given by Kreuzer and Skarke 
\cite[proof of Theorem 1]{KS}, though it requires some work to see that the
condition stated in \cite[Theorem 1]{KS} is equivalent to (C1).

(vii) In theorem \ref{t2.2} (c) conditions $(\N_0^J)_k\neq \emptyset$ for some 
$k\in\N_0$ arise. For $k\in\Z$ denote $\Z^J$, $(\Z^n)_k$ and $(\Z^J)_k$ analogously
to $\N_0^J$, $(\N_0^n)_k$ and $(\N_0^J)_k$. Then $(\Z^J)_k\neq \emptyset$ is equivalent to 
$\gcd(v_j\, |\, j\in J)\, |\, k$. But $(\N_0^J)_k\neq \emptyset$ (for $k\geq 0$) is more 
delicate. In the case $J=\{1,2\}$ sufficient conditions are $\gcd(v_1,v_2)\, |\, k$ and
$\lcm(v_1,v_2)-v_1-v_2+1\leq k$, because then 
$$\left(\frac{\lcm(v_1,v_2)}{v_1}-1\right)\cdot v_1+(-1)\cdot v_2 
= (-1)\cdot v_1+\left(\frac{\lcm(v_1,v_2)}{v_2}-1\right)\cdot v_2$$
is the largest multiple of $\gcd(v_1,v_2)$ missing in $\N_0\cdot v_1+\N_0\cdot v_2$.
\end{remarks}

For any weight system $(v_1,...,v_n,d)\in\N^{n+1}$ with $v_i<d$ define the rational
function
$$\rho_{(v,d)}(t):=\prod_{i=1}^n(t^d-t^{v_i})(t^{v_i}-1)^{-1}.$$
It is well known that $\rho_{(v,d)}(t)\in\N_0[t]$
if a quasihomogeneous polynomial with isolated singularity at $0$ exists.

The conditions $\rho_{(v,d)}(t)\in\Z[t]$ and $\rho_{(v,d)}(t)\in\N_0[t]$ are 
in general weaker than (C1) to (C3), 
but $\rho_{(v,d)}(t)\in\Z[t]$ is equivalent to a surprisingly similar statement.
Denote by $\overline{(C1)}$ and $\overline{(C2)}$
the conditions obtained from (C1) and (C2) in lemma \ref{t2.1} with 
$\N_0$ replaced by $\Z$ in (C1) and (C2) and in the definitions of $R$ and $R_k$.

\begin{lemma}\label{t2.4}
Fix a weight system $(v_1,...,v_n,d)\in\N^{n+1}$ with $v_i<d$.
The following conditions are equivalent.
\begin{list}{}{}
\item[$\overline{(IS3)}$:]
\quad $\rho_{(v,d)}(t)\in\Z[t]$.
\item[(GCD):] 
\quad $\forall \ J\subset N$ the $\gcd(v_j\, |\, j\in J)$ divides \\
\hspace*{4cm} at least $|J|$ of the numbers $d-v_k$.
\item[$\overline{(C2)}$:]
\quad for $R=(\Z^n)_d$. 
\item[$\overline{(C1)}$:]
\quad for $R=(\Z^n)_d$. 
\end{list}
\end{lemma}

{\bf Proof:}
$\overline{(IS3)}$ means that all zeros of $\prod_{i=1}^n(t^{v_i}-1)$ are zeros 
of $\prod_{i=1}^n(t^{d-v_i}-1)$ with at least the same multiplicity. 
This shows $\overline{(IS3)}\iff$(GCD). 
The equivalence (GCD)$\iff\overline{(C2)}$ is trivial.
The equivalence $\overline{(C2)}\iff\overline{(C1)}$ follows as in lemma \ref{t2.1}.
\hfill$\Box$

\begin{lemma}\label{t2.5}
Fix a weight system $(v_1,...,v_n,d)\in\N^{n+1}$ with $v_i<d$.
If $n\leq 3$ then (IS3)$\iff\overline{(IS3)}$.
\end{lemma}

{\bf Proof:}
We restrict to the case $n=3$. It is sufficient to show 
$\overline{(C1)}\Rightarrow$(C1) for $R=(\N_0^n)_d$. 
The $(|J|=1)$-parts of $\overline{(C1)}$ and (C1) coincide.

Consider $J=\{1,2\}$, $J_1=\{1\}$, $J_2=\{2\}$. Then $J$ satisfies (C1) if and only if
$(\N_0^J)_d\neq\emptyset$. 
Now consider the different possibilities how $J_1$ and $J_2$ can satisfy (C1).
The only case where $(\N_0^J)_d\neq\emptyset$ is not obvious 
is the case when $J_1(=\{1\})\& K_1=\{3\}$ and $J_2(=\{2\})\& K_2=\{3\}$ satisfy (C1), 
that is,
when $v_1|(d-v_3)$ and $v_2|(d-v_3)$. Of course, then also $\lcm(v_1,v_2)|(d-v_3)$
and $\lcm(v_1,v_2)\leq d-v_3$.

$\overline{(C1)}$ for $J$ gives $(\Z^J)_d\neq\emptyset$, that is, $\gcd(v_1,v_2)|d$.

The conditions $\gcd(v_1,v_2)|d$ and $\lcm(v_1,v_2)\leq d-v_3$
imply $(\N_0^J)_d\neq\emptyset$ by remark \ref{t2.3} (vii), so $J$ satisfies (C1).

The $(|J|=3)$-part of (C1) follows from the $(|J|=2)$-part.
\hfill$\Box$

\begin{remarks}\label{t2.6}
(i) Lemma \ref{t2.5} is Theorem 3 in \cite{Sa2}. It is also stated in 
\cite[remark after cor. 4.13]{Ar} and \cite[2nd remark in 12.3]{AGV}.

(ii) For $n\geq 4$ $\overline{(IS3)}$ is weaker than (IS3). \cite[12.3]{AGV} contains
the example $(v_1,v_2,v_3,v_4,d)=(1,33,58,24,265)$ of Ivlev.
Here $\rho_{(v,d)}(t)\in\N_0[t]$, but (C1) fails for $J=\{2,4\}$.

(iii) The equivalence $\overline{(IS3)}\iff\overline{(C1)}$ in lemma \ref{t2.4}
is (up to rewriting their condition as $\overline{(C1)}$) Lemma 1 in \cite{KS}.

(iv) Chapter 3 in \cite{Wa} contains results and short proofs for $0$-dimensional
quasihomogeneous complete intersections which are very close to
theorem \ref{t2.2} (b)+(c), lemma \ref{t2.4} and lemma \ref{t2.5}.
\end{remarks}

\section{Types and graphs of quasihomogeneous singularities}\label{s3}
\setcounter{equation}{0}

\noindent
Here a classification of quasihomogeneous polynomials with isolated singularity at $0$
by certain {\it types}, which are encoded in certain {\it graphs}, will be given.
For $n\in\{2,3\}$ this is treated in \cite{Ar}\cite{AGV}, the general case is carried out
in a part of \cite{OR1} which is not published in \cite{OR2}.

The type will come from some choice. Often several choices are possible, and they may
lead to different types or the same type, so, often there are several types
for one quasihomogeneous polynomial.

Now consider $n\in\N$, $N=\{1,...,n\}$, a weight system $(v_1,...,v_n,d)\in \N^{n+1}$
with $v_i<d$ and a quasihomogeneous polynomial $f\in \C[x_1,...,x_n]$
with an isolated singularity at $0$.
Then $\supp f\subset (\N_0^n)_d$ satisfies (C2) by theorem \ref{t2.2} (a).

The choice is a map $\kappa:N\to N$ such that for any $j\in N$ the sets
$J=\{j\}$ and $K=\{\kappa(j)\}$ satisfy (C2) with $R=\supp f$, that is,
$f$ contains a summand $b\cdot x_j^a\cdot x_{\kappa(j)}$ for some $b\in\C^*$, $a\in \N  $.  
The {\it type} is the conjugacy class of this map $\kappa$ with respect to the
symmetric group $S_n$. The graph which encodes the map $\kappa$ is the ordered graph
with $n$ vertices with numbers $1,...,n$ and an arrow from $j$ to $\kappa(j)$ for any
$j\in N$ with $j\neq \kappa(j)$. The ordered graph without the numbering of the
vertices obviously encodes the type.

In order to describe the graphs, an oriented tree is called {\it globally oriented}
if each vertex except one has exactly one outgoing arrow. 
Then the exceptional vertex has only incoming arrows and is called {\it root}.
Starting at any vertex and following the arrows one arrives at the root.

An oriented cycle is called {\it globally oriented} if each vertex has one incoming
and one outgoing arrow. Following the arrows one runs around the cycle. 
The following lemma is obvious.

\begin{lemma}\label{t3.1}
Exactly those graphs occur as graphs of maps $\kappa:N\to N$ whose components 
either are globally oriented trees or consist of one globally oriented cycle and 
finitely many globally oriented trees whose roots are on the cycle.
\end{lemma}

\begin{examples}\label{t3.2}
(i) $n=2$: \cite{Ar}\cite{AGV} 3 types,
\begin{eqnarray*}
\begin{matrix}
\GrBi \ & \GrBii & \GrBiii\\ 
I & II & III 
\end{matrix}
\end{eqnarray*}

(ii) $n=3$: \cite{Ar}\cite{AGV} 7 types.
The sets $J$ under the graphs III and VI are explained in example \ref{t3.6}.
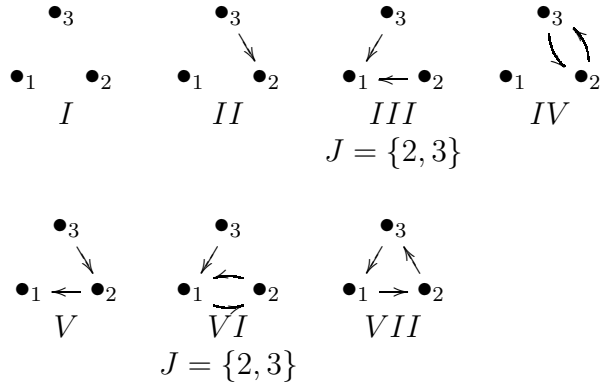
\begin{eqnarray*}
\begin{matrix}
\GrCi \ & \GrCii & \GrCiii & \GrCiv \\ 
I & II & III  & IV \\
 & & J=\{2,3\} & \\ \\
\GrCv & \GrCvi & \GrCvii &  \\ 
V & VI & VII & \\
 & J=\{2,3\} & & 
\end{matrix}
\end{eqnarray*}

(iii) $n=4$: \cite{OR1} and \cite{YS} 19 types.
We follow the numbering in \cite[Proposition 3.5]{YS}.
The sets $J$ under 9 of the 19 graphs are explained in example \ref{t3.6}.

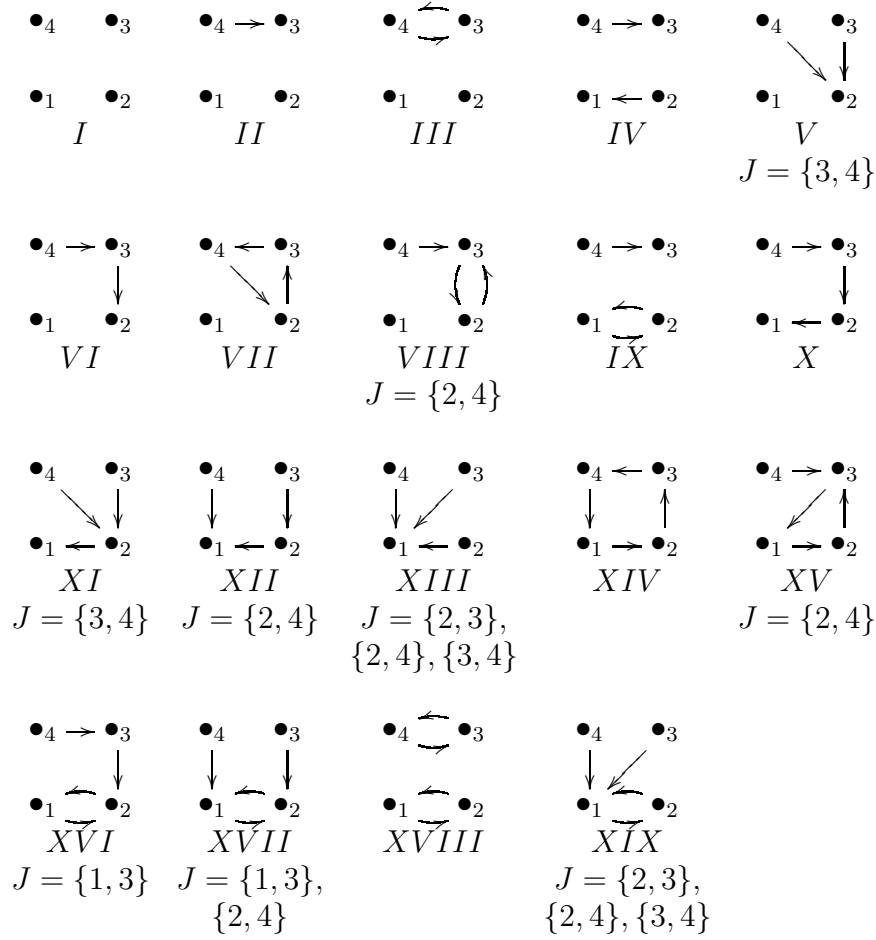
\noindent
\begin{eqnarray*}
\begin{matrix}
\GrDi & \GrDii & \GrDiii & \GrDiv & \GrDv \\
 I & II & III & IV & V \\
 & & & & J=\{3,4\} \\ \\
\GrDvi & \GrDvii & \GrDviii & \GrDix & \GrDx \\
 VI & VII & VIII & IX & X \\
 & & J=\{2,4\} & & \\ \\
\GrDxi & \GrDxii & \GrDxiii & \GrDxiv & \GrDxv  \\
 XI & XII & XIII & XIV & XV \\
 J=\{3,4\} & J=\{2,4\} & J=\{2,3\}, & & J=\{2,4\} \\
 & & \{2,4\}, \{3,4\} & & \\ \\
\GrDxvi & \GrDxvii & \GrDxviii & \GrDxix & \\
 XVI & XVII & XVIII & XIX & \\
 J=\{1,3\} & J=\{1,3\}, & & J=\{2,3\}, &  \\
 & \{2,4\} &  &  \{2,4\}, \{3,4\} &  
\end{matrix}
\end{eqnarray*}

(iv) $n=5$: 47 types.

(v) $n=6$: 128 types.
\end{examples}

\begin{remark}\label{t3.3}
Fix a weight system $(v_1,...,v_n,d)\in \N^{n+1}$ with $v_i<d$, a quasihomogeneous
polynomial $f$ and a map $\kappa:N\to N$ as above. Then for any $j\in N$ the sets
$J=\{j\}$ and $K=\{\kappa(j)\}$ satisfy (C2) with $R=\supp f$ {\it in a unique way}:
There is a unique $a_j\in \N$ with $\alpha:=a_je_j\in R_{\kappa(j)}\cap \N_0^J$,
that is, there is a unique monomial $x_j^{a_j}x_{\kappa(j)}$ 
with exponent $a_je_j+e_{\kappa(j)}$ in the support of $f$.
\end{remark}

Now we forget $(v_1,...,v_n,d)$ and $f$ and start anew with such a tuple of monomials.
We fix $n\in \N$, $N=\{1,...,n\}$, a map $\kappa:N\to N$, numbers $a_1,...,a_n\in \N$
and the set $R:=\{a_je_j+e_{\kappa(j)}\ |\ j\in N\}\subset \N_0^n$
of exponents of the monomials $x_j^{a_j}x_{\kappa(j)}$.

Always $|R|\leq n$, and most often $|R|=n$. The difference $n-|R|$ is the number
of 2-cycles in the graph of $\kappa$ with vertices $j_1$ and $j_2$ and numbers
$a_{j_1}=a_{j_2}=1$.

\begin{lemma}\label{t3.4}
A weight system $(v_1,...,v_n,d)\in \N^{n+1}$ with $v_i<d$ and $R\subset (\N_0^n)_d$
exists if and only if any even cycle with vertices $j_1,...,j_l$ ($l$ even) satisfies
either
\begin{list}{}{}
\item[(EC1)]
neither $a_{j_1}=a_{j_3}=...=a_{j_{l-1}}=1$\\ \hspace*{1cm} nor $a_{j_2}=a_{j_4}=...=a_{j_l}=1$,
\end{list}
or 
\begin{list}{}{}
\item[(EC2)]
$a_{j_1}=a_{j_2}=...=a_{j_l}=1$
\end{list}
(here EC stands for Even Cycle).
If such a weight system exists it is unique up to rescaling if and only if all even
cycles satisfy (EC1).
\end{lemma}

{\bf Proof:}
We work with a normalized weight system $(w_1,...,w_n,1)\in (\Q\cap (0,1))^n\times\{1\}$.
It is a solution of the system of linear equations $a_jw_j+w_{\kappa(j)}=1$, $j\in N$.
We discuss in this order (1) roots of trees not on a cycle,
(2) vertices on a cycle, (3) vertices on trees different from the roots.

\medskip
{\bf (1)} If $j$ is the root of a tree and is not on a cycle then $\kappa(j)=j$ and
$w_j=\frac{1}{a_j+1}\in\Q\cap (0,1)$.

\medskip
{\bf (2)} The restriction of the equations $a_jw_j+w_{\kappa(j)}=1$, $j\in N$, to the vertices
$j_1,...,j_l$ of a cycle with $\kappa(j_l)=j_1$ and $\kappa(j_i)=j_{i+1}$ for $1\leq i\leq l-1$
has a unique solution $(w_{j_1},...,w_{j_l})\in\Q^l$ if and only if
\begin{eqnarray*}
0\neq\det\begin{pmatrix}a_{j_1} & 1 & & \\ & a_{j_2} & & \\ 
 & & \ddots & 1 \\ 1 & & & a_{j_l} \end{pmatrix}
=a_{j_1}\cdot ...\cdot a_{j_l} - (-1)^l,
\end{eqnarray*}
that is, if the cycle is odd or does not satisfy (EC2).

In that case one calculates easily that the solution is 
\begin{eqnarray}\label{3.1}
w_{j_i} = \frac{\rho(a_{j_{i+1}},a_{j_{i+2}},...,a_{j_l},a_{j_1},...,a_{j_{i-1}})}
{a_{j_1}\cdot ...\cdot a_{j_l}-(-1)^l}
\end{eqnarray}
where
\begin{eqnarray}\label{3.2}
\rho:\bigcup_{k=0}^\infty \Z^k &\to & \Z,\\
\rho(x_1,...,x_k)&=&x_1...x_k-x_2...x_k+...+(-1)^{k-1}x_k+(-1)^k.
\end{eqnarray}
If all $x_i\geq 1$ then $\rho(x_1,...,x_k)\geq 0$, and then $\rho(x_1,...,x_k)=0$ 
if and only if $k$ is odd and $x_1=x_3=...=x_k=1$. 
Therefore in the case $a_{j_1}\cdot ...\cdot a_{j_l}-(-1)^l\neq 0$ all $w_{j_i}>0$
if and only if the cycle is odd or is even and satisfies (EC1). 
In that case the inequalities $w_{j_i}>0$ and $w_{j_{i+1}}>0$ and the equation
$a_{j_i}w_{j_i}+w_{j_{i+1}}=1$ show also $0<w_{j_i}$ and $0<w_{j_{i+1}}$.

In the case $a_{j_1}\cdot ...\cdot a_{j_l}-(-1)^l=0$ the cycle is even and satisfies (EC2),
and the equations $a_{j_i}w_{j_i}+w_{j_{i+1}}=1$ give only
\begin{eqnarray*}
w_{j_1}=w_{j_3}=...=w_{j_{l-1}}\quad \textup{and}\quad
w_{j_2}=w_{j_4}=...=w_{j_l}=1-w_{j_1}.
\end{eqnarray*}
Any choice $w_{j_1}\in\Q\cap (0,1)$ works.

\medskip
{\bf (3)} The weights of vertices on the trees different from the roots are successively
determined by
$$w_j=\frac{1-w_{\kappa(j)}}{a_j}$$
and automatically satisfy $0<w_j<1$.
\hfill$\Box$

\bigskip

The following lemma \ref{t3.5} is related to the notion of {\it invertible polynomial}
\cite{ET} and is known to some specialists.
We keep the situation after remark \ref{t3.3}. We need some notations.

The map $\kappa:N\to N$ is of {\it Fermat type} if $\kappa=\id$, that is, if its graph
has no arrows. It is of {\it cycle type} if its graph is a cycle.
It is of {\it chain type} if it has the vertices $j_1,...,j_n$ and the $n-1$ arrows
from $j_i$ to $j_{i+1}$ for $1\leq i\leq n-1$.
The type of $\kappa$ is a {\it sum of Fermat type, cycle types and chain types} if its graph
is a union of the corresponding graphs.

\begin{lemma}\label{t3.5}
Let $n\in \N$, $N=\{1,...,n\}$, $\kappa:N\to N$, $a_1,...,a_n\in\N$ and 
$R=\{a_je_j+e_{\kappa(j)}\, |\, j\in N\}$ be as above such that
any even cycle in the graph of $\kappa$ satisfies (EC1) or (EC2) (in lemma \ref{t3.4})
and such that $a_j\geq 2$ for any $j$ in a tree of the graph of $\kappa$.
Let $(v_1,...,v_n,d)\in \N^{n+1}$ be a weight system with $v_i<d$ 
and $R\subset (\N_0^n)_d$
(it exists by lemma \ref{t3.4}). Then the following 2 conditions are equivalent:
\begin{list}{}{}
\item[(IS4)]
A generic linear combination of the (at most $n$) monomials $x_j^{a_j}x_{\kappa(j)}$,
$j\in N$, is a quasihomogeneous polynomial with an isolated singularity at 0.
\item[(FCC)]
The type of $\kappa$ is a sum of Fermat type, cycle types and chain types.
\end{list}
\end{lemma}

{\bf Proof:}
By theorem \ref{t2.2} (b), (IS4) is equivalent to (C2) for $R$ as above. The implication
(FCC)$\Rightarrow$(IS4) is well known, also a direct proof of (FCC)$\Rightarrow$(C2)
is easy.

The other implication (C2)$\Rightarrow$(FCC) will be proved indirectly: Suppose that 
(FCC) does not hold. Then there are two indices $j_1,j_2\in N$ with $j_1\neq j_2$
and $\kappa(j_1)=\kappa(j_2)$. The set $J:=\{j_1,j_2\}$ does not satisfy (C2) for $R$
as above.\hfill$\Box$

\begin{examples}\label{t3.6}
We return to the examples \ref{t3.2}.

(i) For a fixed map $\kappa:N\to N$ and numbers $a_1,...,a_n\in \N$, the conditions
(EC1) and (EC2) in lemma \ref{t3.4} are not empty if the graph of $\kappa$
contains an even cycle, that is type III for $n=2$, the types IV and VI for $n=3$
and the types III, VIII, IX, XIV, XVI, XVII, XVIII and XIX for $n=4$.

(ii) $n=2$: Type I is Fermat type, type II is chain type, type III is cycle type.
In type III one must avoid $a_1=1,\ a_2>1$ and $a_1>1,\ a_2=1$. Apart from that
(IS4) holds for arbitrary $a_1,a_2\in \N$.

(iii) $n=3$ and $n=4$: 5 of the 7 types with $n=3$ and 10 of the 19 types with $n=4$
are sums of Fermat type, cycle types and chain types. 
There (IS4) and (IS3) (in theorem \ref{t2.2} (c)) hold for
almost arbitrary $a_1,...,a_n\in \N$, with the only constraints
from (EC1) or (EC2) in lemma \ref{t3.4}.

For the other types,
the sets $J$ which fail to satisfy (C1)' for $R=\{a_je_j+e_{\kappa(j)}\ |\ j\in N\}$
are indicated under the graphs in example \ref{t3.2} (ii) and (iii).
For these types one needs more monomials than those with exponents in $R$ 
in order to satisfy (C1)'. This leads to further constraints on the numbers
$a_1,...,a_n$.

(iv) $n=3$: In both cases, III and VI, the failing set is $J=\{2,3\}$.
Suppose that a weight system $(v_1,...,v_n,d)$ as in lemma \ref{t3.4} is determined
from $a_1,a_2,a_3\in \N$ (uniquely except for $a_1=a_2=1$ in type VI).
For (IS3) to hold one needs $(\N_0^J)_d\neq\emptyset$.
By lemma \ref{t2.5} this is equivalent to $(\Z^J)_d\neq\emptyset$ and to
$\gcd(v_1,v_2)\, |\,  d$. This condition is made explicit in \cite{Ar}\cite[13.2]{AGV}.

(v) $n=4$: Suppose that a weight system $(v_1,v_2,v_3,v_4,d)$ as in lemma \ref{t3.4}
is determined from $a_1,a_2,a_3,a_4\in \N$. Consider in each of the 9 cases 
which are not sums of Fermat type, cycle types and chain types a failing set
$J=\{j_1,j_2\}$, that is, with $\kappa(j_1)=\kappa(j_2)=j_3$ and 
$\{1,2,3,4\}=\{j_1,j_2,j_3,j_4\}$. For (IS3) to hold one needs $(\N_0^J)_d\neq\emptyset$
or $(\N_0^J)_{d-v_{j_4}}\neq\emptyset$.
As in lemma \ref{t2.5}, the condition $(\N_0^J)_d\neq\emptyset$ is equivalent to
$(\Z^j)_d\neq\emptyset$ and to $\gcd(v_{j_1},v_{j_2})\, |\, d$. But the condition
$(\N_0^J)_{d-v_{j_4}}\neq\emptyset$ may be stronger than
$(\Z^J)_{d-v_{j_4}}\neq\emptyset$ and $\gcd(v_{j_1},v_{j_2})\, |\, d-v_{j_4}$.

(vi) We consider the case XII with $n=4$ in detail. There 
one starts with arbitrary $a_1,a_2,a_3,a_4\in \N$ and with the monomials
$x_1^{a_1+1}$, $x_2^{a_2}x_1$, $x_3^{a_3}x_2$, $x_4^{a_4}x_1$. The weight system
\begin{eqnarray*}
&&(v_1,v_2,v_3,v_4,d)\\
&&=(a_2a_3a_4,a_1a_3a_4,((a_1+1)(a_2-1)+1)a_4,a_1a_2a_3,(a_1+1)a_2a_3a_4)
\end{eqnarray*}
is unique up to rescaling.
The only failing set is $J=\{2,4\}$, and $\kappa(2)=\kappa(4)=1$, so 
$(\N_0^J)_{d-v_1}\neq\emptyset$.
One needs $(\N_0^J)_d\neq\emptyset$ or $(\N_0^J)_{d-v_3}\neq\emptyset$ for (IS3) to hold.
Now
\begin{eqnarray*}
(\N_0^J)_d\neq\emptyset &\iff& (\Z^J)_d\neq\emptyset \\
&\iff& \gcd(v_2,v_4)\, |\,  d \iff a_1\, |\, \lcm(a_2,a_4).
\end{eqnarray*}
And 
\begin{eqnarray*}
(\N_0^J)_{d-v_3}\neq\emptyset &\Longrightarrow& (\Z^J)_{d-v_3}\neq\emptyset\\
&\iff& \gcd(v_2,v_4)\, |\,  d-v_3\\
&\iff& a_1a_3\, |\,  \frac{a_4}{\gcd(a_2,a_4)}(((a_1+1)(a_2a_3-a_2+1)-1).
\end{eqnarray*}

(vii) Ivlev's example (remark \ref{t2.6} (ii), \cite[12.3]{AGV})
$(v_1,v_2,v_3,v_4,d)=(1,33,58,24,265)$ is of type XII with the monomials
$x_1^{265}$, $x_2^8x_1$, $x_3^4x_2$, $x_4^{11}x_1$, so 
$(a_1,a_2,a_3,a_4)=(264,8,4,11)$. Here $(\N_0^J)_d =\emptyset$ and 
$(\N_0^J)_{d-v_3}=\emptyset$, so (IS3) does not hold, but 
$(\Z^J)_{d-v_3}\neq\emptyset$, so $\oooo{(IS3)}$ holds and $\rho_{{\bf v},d}(t)\in\Z[t]$,
even $\in\N_0[t]$.
\end{examples}

Two function germs $f_1,f_2\in\OO_{\C^n,0}$ are right equivalent if there is a
local coordinate change $\varphi:(\C^n,0)\to(\C^n,0)$ such that $f_1\circ \varphi=f_2$.
Often in one right equivalence class of functions with an isolated singularity at $0$,
there are several quasihomogeneous functions with different weight systems.
For example $x_1^{a_1}x_2+x_2^{a_2}x_3+x_3x_1$ with weight system
$(v_1,v_2,v_3,d)=(a_2,1,a_1a_2-a_2+1,a_1a_2+1)$ and $x_1^{a_1a_2+1}+x_2^2+x_3^2$
with weight system $(v_1',v_2',v_3',v_3',d')=(2,a_1a_2+1,a_1a_2+1,2a_1a_2+2)$ are in the 
same right equivalence class of $A_{a_1a_2}$-singularities \cite{ET}.
The ambiguity was analysed in \cite{Sa1}.

\begin{theorem}\label{t3.7}\cite{Sa1}
Let $f\in\OO_{\C^n,0}$ be a function germ with an isolated singularity at 0.

(a) $f$ is right equivalent to a quasihomogeneous polynomial if and only if
$$f\in J_f:=\left(\frac{\paa f}{\paa x_1},...,
\frac{\paa f}{\paa x_n}\right)\subset \OO_{\C^n,0}.$$

(b) If $f$ is quasihomogeneous with normalized weight system $(w_1,...,w_n,1)$ with 
$0<w_1\leq ...\leq w_n<1$ and if $f\in{\bf m}^3_{\C^n,0}$, then the weight system
is unique and $0<w_1\leq ...\leq w_n<\frac{1}{2}$.

(c) If $f\in J_f$ then $f$ is right equivalent to a quasihomogeneous polynomial
$g(x_1,...,x_k)+x_{k+1}^2+...+x_n^2$ with $g\in{\bf m}^3_{\C^k,0}$.
Especially, its normalized weight system satisfies 
$0<w_1\leq ...\leq w_k<w_{k+1}=...=w_n=\frac{1}{2}$.

(d) If $f$ and $\www f\in \OO_{\C^n,0}$ are right equivalent and quasihomogeneous
with normalized weight systems 
$(w_1,...,w_n,1)$ and $(\www w_1,...,\www w_n,1)$ with 
$w_1\leq ...\leq w_n\leq\frac{1}{2}$ and $\www w_1\leq ...\leq \www w_n\leq\frac{1}{2}$
then $w_i=\www w_i$.
\end{theorem}

\begin{remarks}\label{t3.8}
(i) Part (b) can be proved with the arguments in the proof of lemma \ref{t3.4}.
The condition $f\in{\bf m}^3_{\C^n,0}$ is equivalent to all $a_2\geq 2$.
The implication $w_j<\frac{1}{2}$ is nontrivial only in case (2) in the proof
of lemma \ref{t3.4}.

(ii) Part (c) follows from (a) and the splitting lemma and (b).

(iii) An argument for part (d) different from the proof in \cite{Sa1} is as follows.
If $f$ is quasihomogeneous with some weight system $(v_1,...,v_n,d)\in\N^{n+1}$
then $\rho_{({\bf v},d)}(t)\in\N_0[t]$, so
$$\rho_{({\bf v},d)}(t^{1/d})=\sum_{j=1}^\mu t^{\alpha_j}$$
for certain numbers $\alpha_1,...,\alpha_\mu\in\frac{1}{d}\N$.
These numbers and $\rho_{({\bf v},d)}(t^{1/d})$ are invariants of the right equivalence
class of $f$. This is well known and follows essentially from calculations in \cite{Br}.
The numbers $\alpha_1,...,\alpha_\mu$ are the {\it exponents} of the right equivalence
class of $f$.
By part (c) there exists a weight system $(\www v_1,...\www v_n,\www d)$ 
with $\www v_i\leq\frac{\www d}{2}$ and 
$$\sum_{j=1}^\mu t^{\alpha_j} =\rho_{({\bf \www v},\www d)}(t^{1/\www d}).$$
It is easy to see that one can recover the normalized weight system 
$\frac{1}{\www d}(\www v_1,...,\www v_n,\www d)$ from the exponents and this 
equation. Therefore this normalized weight system is unique.
\end{remarks}

\section{Milnor number versus weighted degree}\label{s4}
\setcounter{equation}{0}

\noindent
Let $p_i$, $i\in\N$, be the $i$-th prime number, so $(p_1,p_2)=(2,3)$.
Define 
$$l(n):=\prod_{i=1}^n\frac{p_i}{p_i-1},$$
so $(l(1),l(2),l(3),l(4),l(5))=(2,3,\frac{15}{4},\frac{35}{8},\frac{77}{16})$.
The prime number theorem in the form $p_n=n\log n\cdot (1+o(1))$ \cite[Theorem 8]{HW}
and Mertens' theorem 
$$\prod_{prime\ numbers\ p\leq x}\frac{p}{p-1}=e^\gamma\cdot\log x\cdot (1+o(1))$$
with $\gamma=$ Euler's constant \cite[Theorem 429]{HW} imply
$$l(n)=e^\gamma\cdot\log n\cdot (1+o(1)).$$

\begin{theorem}\label{t4.1}
(a) Let $f\in\C[x_1,...,x_n]$ be a quasihomogeneous polynomial with 
an isolated singularity at $0$ and 
reduced weight system $(v_1,...,v_n,d)\in\N^{n+1}$ with $v_i\leq\frac{d}{2}$
for all $i$ (reduced: $\gcd(v_1,...,v_n,d)=1$). Then
$$d\leq l(n)\cdot\mu.$$
(b) If $v_i<\frac{d}{2}$ for all $i$ and $n\geq 2$ then
$$d\leq l(n-1)\cdot \mu.$$
\end{theorem}

These estimates rely only on the conditions for $J$ with $|J|=1$ in (C1)-(C3)
for $R=(\N_0^n)_d$, the conditions for $|J|\geq 2$ are not needed.
Theorem \ref{t4.3} formulates this more general case.
Both theorems are proved after stating theorem \ref{t4.3}.

\begin{remarks}\label{t4.2}
(i) These estimates are useful for a classification of such weight systems
using computer, for a fixed number of variables and with Milnor numbers 
up to a chosen bound. See section \ref{s5}.

(ii) Calculations in \cite{Br} show that for a quasihomogeneous polynomial $f$ as
in theorem \ref{t4.1} the monodromy on the Milnor lattice is semisimple
with eigenvalues $e^{-2\pi i\alpha_1}$, ..., $e^{-2\pi i\alpha_\mu}$,
where $\alpha_1$,..., $\alpha_\mu$ are the exponents considered in remark \ref{t3.8} (iii).
For $f\in{\bf m}^3_{\C^n,0}$ the procedure mentioned in remark \ref{t3.8} (iii), 
which recovers the normalized weights $(w_1,...,w_n)$ from the exponents, shows 
that the tuples $(w_1,...,w_n)$ and $(\alpha_1,...,\alpha_\mu)$ have the same
common denominator $d$. Therefore in the case $f\in{\bf m}^3_{\C^n,0}$ the 
order of the monodromy is $d$. 
Adding squares $x_{n+1}^2+...x_{n+m}^2$ changes the eigenvalues by the factor
$(-1)^m$ and replaces $d$ by $\www d$ with $\www d=2d$ for odd $d$ and $\www d=d$
for even $d$. Then the order of the monodromy is $\www d$ or $\frac{\www d}{2}$.
\end{remarks}

\begin{theorem}\label{t4.3}
Fix $n\in \N$, $N=\{1,...,n\}$, a map $\kappa:N\to N$, numbers $a_1,...,a_n\in\N$
and the set $R=\{a_je_j+e_{\kappa(j)}\, |\, j\in N\}$ of exponents of the 
monomials $x_j^{a_j}x_{\kappa(j)}$.
Suppose that $a_j\geq 2$ for all $j\in N$ which lie in components $C$ of the
graph of $\kappa$ with $|C|\geq 2$.

(a) By lemma \ref{t3.4} there is a unique reduced weight system $(v_1,...,v_n,d)\in\N^{n+1}$
with $R\subset (\N_0^n)_d$. It satisfies $v_j<\frac{d}{2}$ for $a_j\geq 2$ and 
$v_j=\frac{d}{2}$ for $a_j=1$. Define
$$\mu:=\prod_{j=1}^n\left(\frac{d}{v_j}-1\right).$$

(b) $$d\leq l(n)\cdot \mu.$$

(c) If all $a_j\geq 2$ and $n\geq 2$ then 
$$d\leq l(n-1)\cdot \mu.$$

(d) If $n=1$ then $d=a_1+1$ and $\mu=a_1$.
\end{theorem}

{\bf Proof of theorem \ref{t4.1}:}
Suppose $v_1\leq ...\leq v_k<v_{k+1}=...=v_n=\frac{1}{2}$ for some $k$ with $0\leq k\leq n$.
By theorem \ref{t3.7} $f$ is right equivalent to a quasihomogeneous polynomial
$g(x_1,...,x_k)+x_{k+1}^2+...+x_n^2$ with $g\in{\bf m}^3_{\C^k,0}$ with an 
isolated singularity at $0$ and the same weight system $(v_1,...,v_n,d)$.

Choose a map $\kappa:N\to N$ for $g+x_{k+1}^2+...+x_n^2$ as in section \ref{s3}.
By remark \ref{t3.3} there are unique numbers $a_1,...,a_n\in\N$ such that
$a_je_j+e_{\kappa(j)}$ are in $\supp(g+x_{k+1}^2+...+x_n^2)$.
The hypotheses in theorem \ref{t4.3} are satisfied. Theorem \ref{t4.3} (b) and (c)
give theorem \ref{t4.1} (a) and (b).
\hfill$\Box$

\bigskip
{\bf Proof of theorem \ref{t4.3}:}
(a) The first part follows from lemma \ref{t3.4}. If $a_j=1$ then $j$ is itself a component
of the graph of $\kappa$, so $(a_j+1)v_j=d$, so $v_j=\frac{d}{2}$.
If $a_j\geq 2$ then $j$ lies in a component $C$ of the graph of $\kappa$ with $a_i\geq 2$
for all $i\in C$. Then $v_j<\frac{d}{2}$ follows as in remark \ref{t3.8} (i) with the
arguments in the proof of lemma \ref{t3.4}.

\medskip

(b) and (c)
Write $\frac{v_j}{d}=w_j=\frac{s_j}{t_j}$ with $w_j\in\Q\cap (0,\frac{1}{2}]$ and
$s_j,t_j\in\N$, $\gcd(s_j,t_j)=1$. An elementary, but important observation is 
\begin{eqnarray}\label{4.1}
j\neq\kappa(j)\Longrightarrow t_j=t_{\kappa(j)}\cdot \beta_j
\textup{ for some }\beta_j\in \N\textup{ with }\beta_j\, |\, a_j.
\end{eqnarray}
This follows from
\begin{eqnarray*}
\frac{s_j}{t_j}=w_j=\frac{1-w_{\kappa(j)}}{a_j}
=\frac{t_{\kappa(j)}-s_{\kappa(j)}}{t_{\kappa(j)}\cdot a_j}
\quad\textup{and}\quad \gcd(t_{\kappa(j)},t_{\kappa(j)}-s_{\kappa(j)})=1.
\end{eqnarray*}

For any subset $C\subset N$ define
\begin{eqnarray*}
\mu(C)&:=&\prod_{j\in C}(\frac{1}{w_j}-1),\quad \textup{especially }
\mu(\emptyset)=1,\ \mu(N)=\mu,\\
d(C)&:=&\lcm(t_j\, |\, j\in C),\quad \textup{especially }
d(\emptyset)=1,\ d(N)=d.
\end{eqnarray*}

Let $C_{Fermat}$ be the union of all components $C$ of the graph of $\kappa$ with $|C|=1$.
For $j\in C_{Fermat}$ $w_j=\frac{1}{a_j+1}$, so
\begin{eqnarray}\label{4.2}
\mu(C_{Fermat}) &=& \prod_{j\in C_{Fermat}}a_j,\\
d(C_{Fermat}) &=& \lcm(a_j+1 \, |\, j\in C_{Fermat}).\label{4.3}
\end{eqnarray}

Now we will study $\mu(C)$ and $d(C)$ for a component $C$ of the
graph of $\kappa$ with $|C|\geq 2$. By hypothesis $a_j\geq 2$ for $j\in C$.

\medskip
{\bf Case 1,} $C$ is a cycle:
Suppose $C=\{1,...,m\}$ with $\kappa(j)=j-1$ for $2\leq j\leq m$ and $\kappa(1)=m$.
\eqref{4.1} gives immediately $t_1=t_2=...=t_m=d(C)$.
\eqref{3.1} shows (with $\rho $ as in \eqref{3.2})
\begin{eqnarray}\label{4.4}
d(C) &=& t_1=...=t_m=\frac{1}{\gamma}\cdot (a_1...a_m-(-1)^m)\\
\textup{where }
\gamma &=& \gcd (a_1...a_m -(-1)^m, \rho(a_{j-1},...,a_1,a_m,...,a_{j+1}))\label{4.5}
\end{eqnarray}
for any $j\in\{1,...,m\}$. 
Define here $\www d(C):=\gamma\cdot d(C)=a_1...a_m-(-1)^m$. 

One calculates
\begin{eqnarray}\nonumber
\mu(C)&=& \prod_{j=1}^m\frac{d-v_j}{v_j} 
=\prod_{j=1}^m \frac{a_1...a_m-(-1)^m-\rho(a_{j-1},...,a_1,a_m,...a_{j+1})}
{\rho(a_{j-1},...,a_1,a_m,...,a_{j+1})}\\
&=&\prod_{j=1}^m\frac{a_{j+1}\cdot \rho(a_j,...,a_1,a_m,...,a_{j+2})}
{\rho(a_{j-1},...,a_1,a_m,...,a_{j+1})}
=a_1\cdot ...\cdot a_m.\label{4.6}
\end{eqnarray}

\medskip
{\bf Case 2,} $C$ is not a cycle: 
Then $C$ is either a tree or a cycle with one or several attached trees.
If $C$ is a tree suppose $C_1=\{1\}\subset C$ is the root, and define $m:=1$.
If $C$ is a cycle with attached trees suppose $C_1=\{1,...,m\}$ is the cycle,
and $\kappa(j)=j-1$ for $ 2\leq j\leq m$, $\kappa(1)=m$.
In both cases the {\it set of leaves} is the subset $C_2\subset C-C_1$ of vertices
with no incoming arrows. For any leaf $j\in C_2$ denote by $C(j)$ the set of vertices
on the path from $j$ to $C_1$, excluding the vertex in $C_1$, so
\begin{eqnarray*}
C(j)=\{j,\kappa(j),...,\kappa^{l(j)}(j)\}\subset C-C_1
\textup{ with }\kappa^{l(j)+1}(j)\in C_1.
\end{eqnarray*}
Then with $\gamma:=1$ if $m=1$ and $\gamma$ as in \eqref{4.5} if $m\geq 2$ one has
\begin{eqnarray*}
d(C_1)&=&\frac{1}{\gamma}\cdot (a_1...a_m-(-1)^m).
\end{eqnarray*}
With \eqref{4.1} and $\beta_i$ as defined in \eqref{4.1} one finds
\begin{eqnarray}
t_j&=& d(C_1)\cdot \prod_{i\in C(j)}\beta_i\qquad \textup{for }j\in C_2,\label{4.7}\\
d(C)&=& \lcm(t_j\, |\, t_j\in C_2)\nonumber\\
&=& d(C_1)\cdot \lcm(\prod_{i\in C(j)}\beta_i\, |\,  j\in C_2).\label{4.8}
\end{eqnarray}
We will estimate $d(C)$ by $\www d(C)$ with $d(C)\, |\, \www d(C)$ and 
\begin{eqnarray}
\www d(C) := (a_1...a_m-(-1)^m)\cdot\left(\prod_{j\in C-(C_1\cup C_2)}a_j\right)\cdot
\lcm(a_j\, |\, j\in C_2).\label{4.9}
\end{eqnarray}

In order to estimate $\mu(C)$ from above, we choose a decomposition of $C-C_1$ into a
disjoint union
\begin{eqnarray*}
C-C_1 =\bigcup_{j\in C_2}^{.}\www C(j)
\end{eqnarray*}
with $\www C(j)\subset C(j)$ being a suitable sub-chain of $C(j)$,
\begin{eqnarray*}
\www C(j) =\{j,\kappa(j),...,\kappa^{\www l(j)}(j)\}
\quad\textup{ for some }\www l(j)\leq l(j).
\end{eqnarray*}
To simplify notations suppose for a moment that one such sub-chain $\www C(j)$ takes the
form $\www C(j)=\{j,j-1,...,k\}$ with $\kappa(i)=i-1$ for $k\leq i\leq j$. 
Using $w_l=\frac{1-w_{\kappa(l)}}{a_l}$ repeatedly one finds by an easy induction
for $k\leq i\leq j$
\begin{eqnarray}\label{4.10}
w_i =\frac{\rho(a_{i-1},...,a_{k+1},a_k)+(-1)^{i-1-k}w_{k-1}}
{a_k a_{k+1}...a_{i-1} a_i}.
\end{eqnarray}
Therefore 
\begin{eqnarray}\nonumber
\mu(\www C(j))&=& \prod_{i\in \www C(j)}\frac{1-w_i}{w_i}
=\prod_{i\in \www C(j)}\frac{\rho(a_{i},...,a_{k+1},a_k)+(-1)^{i-k}w_{k-1}}
{\rho(a_{i-1},...,a_{k+1},a_k)+(-1)^{i-1-k}w_{k-1}}\\
&=& \frac{\rho(a_{j},...,a_{k+1},a_k)+(-1)^{j-k}w_{k-1}}
{1-w_{k-1}}.\label{4.11}
\end{eqnarray}
Because all $a_i\geq 2$ for $i\in C$, one can estimate
\begin{eqnarray}\nonumber
\rho(a_{j},...,a_{k+1},a_k)+(-1)^{j-k}w_{k-1}>a_k...a_{j-1}\cdot(a_j-1),\\
\mu(\www C(j)) >\frac{a_k... a_{j-1}\cdot(a_j-1)}{1-w_{k-1}}
> a_k... a_{j-1}\cdot(a_j-1).\label{4.12}
\end{eqnarray}

The following additional estimate is relevant only for odd $m$. But it holds for all $m$,
and it will be smoother to treat even and odd $m$ simultaneously. For $k-1\in C_1$
\begin{eqnarray}\nonumber
&&\mu(C_1)\cdot \frac{1}{1-w_{k-1}} \\
&=& a_1...a_m\cdot
\frac{a_1...a_m-(-1)^m}{a_1...a_m-(-1)^m-\rho(a_{k-2},...,a_1,a_m,...,a_k)}\nonumber \\
&\geq& a_1...a_m-(-1)^m.\label{4.13}
\end{eqnarray}

Now we put together the pieces and estimate $\mu(C)$ from above.
There is (at least) one leaf $j_0\in C_2$ with $\www C(j_0)=C(j_0)$,
so $k-1:=\kappa^{\www l(j_0)+1}(j)\in C_1$. For this leaf $j_0$ we use the finer estimate
in \eqref{4.12}
\begin{eqnarray*}
\mu(\www C(j_0))> \frac{1}{1-w_{k-1}}\cdot (a_{j_0}-1)\cdot \prod_{i\in C(j_0)-\{j_0\}}a_i.
\end{eqnarray*}
Together with \eqref{4.12} for all other leaves $j\in C_2$
and \eqref{4.13} we obtain
\begin{eqnarray}\nonumber
&&\mu(C)= \mu(C_1)\cdot \prod_{j\in C_2}\mu(\www C(j))\\
&\geq& (a_1...a_m-(-1)^m)\cdot \left(\prod_{j\in C-(C_1\cup C_2)}a_j\right)
\cdot \left(\prod_{j\in C_2}(a_j-1)\right).\label{4.14}
\end{eqnarray}

\medskip
Now case 2 is finished. We can estimate $d$ and $\mu$ and their quotient.
$C_{Leaf}\subset N-C_{Fermat}$ denotes the union of the leaves of all components
$C$ with $|C|\geq 2$. For any such $C$ the notations of case 2 are preserved,
$C_1$ is the root or the cycle in it, and $C_2$ is the set of leaves in it.
If $C$ is a cycle then $C=C_1$.
\begin{eqnarray}
d &=& \lcm\left(d(C_{Fermat}); \ d(C) 
\textup{ for all components }C\textup{ with }|C|\geq 2\right)
\nonumber\\
&\leq & \lcm\left(d(C_{Fermat}); \ \www d(C) \textup{ for }C\textup{ with }|C|\geq 2
\right)\nonumber\\
&\leq & \prod_{C\textup{ with }|C|\geq 2,\, C\textup{ not an odd cycle}}
\left( \left(\prod_{j\in C_1}a_j-(-1)^{|C_1|}\right)\cdot 
\prod_{j\in C-(C_1\cup C_2)}a_j\right)\cdot \nonumber\\
&& \lcm\Bigl(a_j+1\textup{ for }j\in C_{Fermat}; \nonumber\\
&& \prod_{j\in C}a_j+1\textup{ for }C
\textup{ an odd cycle }; \ a_j\textup{ for }j\in C_{Leaf}\Bigr) .
\label{4.15}
\end{eqnarray}

\begin{eqnarray}\nonumber
\mu&=& \mu(C_{Fermat})\cdot \prod_{C\textup{ a cycle}}\mu(C)\cdot 
\prod_{C\textup{ not a cycle,}|C|\geq 2}\mu(C)\\
&\geq& \prod_{j\in C_{Fermat}}a_j\cdot 
\prod_{C\textup{ a cycle}}\left(\prod_{j\in C}a_j\right)\cdot \label{4.16}\\
&&\prod_{C\textup{ not a cycle,}|C|\geq 2}\left( \left(\prod_{j\in C_1}a_j-(-1)^{|C_1|}\right)
\cdot \prod_{j\in C-(C_1\cup C_2)}a_j\cdot \prod_{j\in C_2}(a_j-1)\right).
\nonumber
\end{eqnarray}

\begin{eqnarray}
\frac{d}{\mu} \leq 
\frac{\lcm\begin{pmatrix}
a_j+1\textup{ for }j\in C_{Fermat}; \\
\prod_{j\in C}a_j+1\textup{ for }C
\textup{ an odd cycle}; \ a_j\textup{ for }j\in C_{Leaf}\end{pmatrix}}
{\prod_{j\in C_{Fermat}}a_j\cdot
\prod_{C\textup{ an odd cycle}}\left(\prod_{j\in C}a_j\right)\cdot 
\prod_{j\in C_{Leaf}}(a_j-1).}
\label{4.17}
\end{eqnarray}
In lemma \ref{t4.4} two numbers $l_1(n)$ and $l_2(n)\in \Q_{>0}$ are defined.
Obviously $\frac{d}{\mu}\leq l_1(n)$, and if all $a_j\geq 2$ and $n\geq 2$ then
$\frac{d}{\mu}\leq \max(l_2(n),l_1(n-1))$.
The parts (b) and (c) of theorem \ref{t4.3} follow now with lemma \ref{t4.4}.
Part (d) is trivial.
\hfill$\Box$

\begin{lemma}\label{t4.4}
For $n\in \N$ define
\begin{eqnarray*}
l_1(n)&=& \max\left( \frac{\lcm(b_1,...,b_n)}{(b_1-1)\cdot ...\cdot (b_n-1)}\ |\ 
b_1,...,b_n\in \N-\{1\}\right),\\
l_2(n)&=& \max\left( \frac{\lcm(b_1,...,b_n)}{(b_1-1)\cdot ...\cdot (b_n-1)}\ |\ 
b_1,...,b_n\in \N-\{1,2\}\right).
\end{eqnarray*}
Then 
\begin{eqnarray*}
l_1(n)=l(n):=\prod_{i=1}^n\frac{p_i}{p_i-1}\geq l_2(n+1),
\end{eqnarray*}
here $p_i$ is the $i$-th prime number.
\end{lemma}

{\bf Proof:}
First, $l_1(n)=l(n)$ will be proved. Choose $b_1,...,b_n\in \N$ arbitrarily.
Write $\lcm(b_1,...,b_n)=\prod_{i\in I}p_i^{r_i}$ with $I\subset \N$ finite,
$r_i\geq 1$ for $i\in I$. For any $i\in I$ choose $\beta(i)\in N$ with
$p_i^{r_i}\, |\, b_{\beta(i)}$. Define
\begin{eqnarray*}
\www b_j:= \prod_{i\textup{ with }\beta(i)=j} p_i^{r_i}.
\end{eqnarray*}
For any $j$ with $\www b_j>1$ let $i(j)$ be the minimal $i$ with $\beta(i)=j$.
Then
\begin{eqnarray*}
\lcm(b_1,...,b_n) &=& \lcm(\www b_j\, |\, \www b_j>1) 
= \prod_{j\textup{ with }\www b_j>1}\www b_j,\\
\frac{\lcm (b_1,...,b_n)}{(b_1-1)\cdot ...\cdot (b_n-1)}
&\leq& \prod_{j\textup{ with }\www b_j>1}\frac{\www b_j}{\www b_j-1}\\
&\leq& \prod_{j\textup{ with }\www b_j>1}\frac{p_{i(j)}}{p_{i(j)}-1}
\leq \prod_{i=1}^n\frac{p_i}{p_i-1}.
\end{eqnarray*}
This proves $l_1(n)\leq l(n)$. The choice $b_i=p_i$ proves $l_1(n)\geq l(n)$.

Analogously one shows for $n\geq 2$
\begin{eqnarray*}
l_2(n) = \frac{3}{3-1}\cdot \frac{4}{4-1}\cdot \prod_{i=3}^n\frac{p_i}{p_i-1}.
\end{eqnarray*}
$l_2(2)=2=l(1)$. For $n\geq 2$ the estimate $l_2(n+1)\leq l(n)$ follows from
$$\frac{4}{4-1}\cdot \frac{p_{n+1}}{p_{n+1}-1}\leq \frac{4}{3}\cdot\frac{p_3}{p_3-1}
=\frac{5}{3}<\frac{2}{2-1}.
$$
\hfill$\Box$

\section{Computer calculations}\label{s5}
\setcounter{equation}{0}

\noindent
Theorem \ref{t2.2} (c) gives combinatorial characterizations (C1)-(C3) of those reduced
weight systems $(v_1,...,v_n,d)\in \N^{n+1}$ for which quasihomogeneous polynomials
with an isolated singularity at $0$ exist. These characterizations can be used in
computer programs to find all such weight systems with Milnor number up to some chosen bound.
Because of theorem \ref{t3.7} for most purposes it is sufficient to restrict to weight
systems with $v_i< \frac{d}{2}$. Theorem \ref{t4.1} (b) gives then the
bound $d\leq l(n-1)\cdot \mu$ for $d$ if $n\geq 2$.

The second author carried out such computer calculations for $n=2,3,4$.
The following table lists for $n=2,3,4$ the number of reduced weight systems
$(v_1,...,v_n,d)$ (up to reordering of $v_1,...,v_n$) with $v_i<\frac{d}{2}$
which satisfy (C1)-(C3) for $R=(\N_0^n)_d$ and whose Milnor number is less or equal than 
the number $\mu$ in the left column. 

\medskip

\begin{tabular}{r|r|r|r|r}
$\mu$ & $n=1$ & $n=2$ & $n=3$ & $n=4$\\ \hline 
50 & 50 & 187 & 217 & 100 \\
100 & 100 & 493 & 806 & 590 \\
150 & 150 & 847 & 1627 & 1442 \\
200 & 200 & 1242 & 2623 & 2678 \\
300 & 300 & 2083 & 5027 & 6059 \\
400 & 400 & 2998 & 7832 & 10459 \\
500 & 500 & 3957 & 10931 & 15634 \\
1000 & 1000& 9246 & 30241 & 52761 \\
1500 & 1500 & 15058 & 53698 & 103841 \\
2000 & 2000 & 21194 & 80055 & 165624 \\
3000 & 3000 & 34177 & 139343 & ? \\
4000 & 4000 & 47833 & 205191 & ? \\
5000 & 5000 & 62012 & 276169 & ? \\
6000 & 6000 & 76545 & 351335 & ? \\
7000 & 7000 & 91439 & 430009 & ? \\
8000 & 8000 & 106616 & 512141 & ? \\
9000 & 9000 & 122040 & 596879 & ? 
\end{tabular}

\medskip

On the homepage \cite{HK} 
tables with all these weight systems and the characteristic polynomials
of the monodromy are available. Of course for $n=1$ one has just the 
$A_\mu$-singularities $x_1^{\mu+1}$ with $(v_1,d)=(1,\mu+1)$ for $\mu\geq 1$. 
The $A_1$-singularity is taken into account in the column for $n=1$ 
despite $v_1=\frac{d}{2}$ in that case.

For example, the total number of reduced weight systems for $n=4$ 
with $v_i\leq \frac{d}{2}$ and (C1)-(C3) and $\mu\leq 50$ is $50+187+217+100.$

The weight system $(\frac{v_1}{d},\frac{v_2}{d},\frac{v_3}{d},\frac{v_4}{d})$
with $\frac{v_i}{d}<\frac{1}{2}$ and the largest $d$ within $\mu\leq 500$ 
is $(\frac{1}{58},\frac{1}{5},\frac{1}{3},\frac{57}{116})$ with 
$\mu=473$, $d=1740$, $l(3)\cdot \mu=1773,75$. This indicates that the estimate in
theorem \ref{t4.1} (b) cannot be improved much.

For any $n$ the weight system with $v_i<\frac{d}{2}$ with the smallest Milnor number
is $(1,...,1,3)$ with $d=3$ and $\mu=2^n$. This follows from \cite[Lemma 2]{KS}.
This lemma says that there is an injective map
$$\nu:\{i\, |\, v_i>\frac{1}{3}\}\to \{i\, |\, v_i<\frac{1}{3}\}
\quad \textup{with}\quad v_{\nu(i)}=d-2v_i.$$
Then
$$(\frac{d}{v_i}-1)(\frac{d}{v_{\nu(i)}}-1) >4.$$

For $n=2$ weight systems with $v_i<\frac{d}{2}$ exist for any $\mu\geq 4$, 
because of the $D_\mu$-singularities $x_1^{\mu-1}+x_2^2x_1$.
But for $n=3$ and $n=4$ there are some gaps, some numbers $>2^n$ which are
not Milnor numbers of any quasihomogeneous singularities $f\in {\bf m}^3_{\C^n,0}$.
We list all gaps up to $\mu=1000$ for $n=3$ and up to $\mu=500$ for $n=4$.
\begin{eqnarray*}
n=3:& \mu=& 9,13,37,61,73, 157, 193, 277, 313, 397, 421,\\
&& 457, 541, 613, 661, 673, 733, 757, 877, 997.\\
n=4:& \mu =& 17,18,19,23,27,47,59,74,83, 107, 167, 179,\\
&& 219, 227, 263, 314, 347, 359, 383, 467, 479.
\end{eqnarray*}
Corollary \ref{t6.3} will give an explanation of the majority
of these gaps in terms of Sophie Germain prime numbers and
similar prime numbers.

Yonemura \cite{Yo} had classified all reduced weight systems 
$(v_1,v_2,v_3,v_4,d)$ with $\sum_i v_i=d$ and (C1)-(C3) for 
$R=(\N_0^n)_d$. Using our lists, we recovered his 95 weight systems.
48 are in our list for $n=3$ with $\sum_{i=1}^3v_i=\frac{d}{2}$,
with Milnor numbers ranging between $125$ ($(1,1,1,6)$ and 
$492$ ($(1,6,14,42)$). 47 are in our list for $n=4$, with
Milnor numbers ranging between $81$ ($(1,1,1,1,4)$) and 
264 ($(1,3,7,10,21)$).

\section{The case Milnor number = prime number}\label{s6}
\setcounter{equation}{0}

\noindent
The computer calculations mentioned in section \ref{s5} led us to expect the following result.
This section is devoted to its proof.

\begin{theorem}\label{t6.1}
Let $f\in\C[x_1,...,x_n]$ be a quasihomogeneous polynomial with an isolated singularity
at $0$ and normalized weight system $(w_1,...,w_n,1)\in(\Q\cap(0,\frac{1}{2}))^n\times\{1\}$
such that its Milnor number $\mu$ is a prime number.

(a) There are numbers $a_1,...,a_n\in\N-\{1\}$ and $c_1,...,c_n\in\C^*$ such that
\begin{eqnarray*}
f=c_1x_1^{a_1+1} + c_2x_2^{a_2}x_1 + ... + c_nx_n^{a_n}x_{n-1}.
\end{eqnarray*}
Therefore $f$ is of chain type by the map $\kappa:N\to N$ with $\kappa(1)=1$,
$\kappa(j)=j-1$ for $2\leq j\leq n$. And this is the only possible map $\kappa$ as
in section \ref{s3}. Also, by rescaling of $x_1,...,x_n$ one can arrange
$c_1=...=c_n=1$. So, $f$ is unique up to right equivalence.

\medskip
(b) Write $w_i=\frac{s_i}{t_i}$ with $s_i,t_i\in\N$, $\gcd(s_i,t_i)=1$.
Then 
\begin{eqnarray*}
t_i&=&a_i...a_2\cdot (a_1+1),\quad d=t_n,\\
s_i&=&\rho(a_{i-1},...,a_2,a_1+1)\quad (\textup{with }\rho\textup{ as in }\eqref{3.2}),\\
s_1&=&1,\quad s_{i+1}=t_i-s_i=t_i-t_{i-1}+t_{i-2}-...+(-1)^i,\\
\mu&=&\rho(a_n,...,a_2,a_1+1).
\end{eqnarray*}

(c) The characteristic polynomial of the monodromy on the Milnor lattice of $f$
is $\prod_{m:\eqref{6.1}}\Phi_m$, here $\Phi_m$ is the cyclotomic polynomial of the 
$m$-th primitive unit roots, and \eqref{6.1} is the condition
\begin{eqnarray}\label{6.1}
m\, |\, a_n...a_2(a_1+1),\quad \min(i\ |\ m\, |\, a_i...a_2(a_1+1))\equiv n\mod 2.
\end{eqnarray}
Especially, all eigenvalues have multiplicity $1$.
\end{theorem}

\begin{examples}\label{t6.2}
For $n=2,3$ all tuples $(a_1,...,a_n)$ as in theorem \ref{t6.1} with $\mu\leq 23$
are listed below, for $n=4$ all tuples with $\mu\leq 31$.

\medskip

\noindent
\begin{tabular}{l|l|l|l}
$\mu$ & $n=2$ & $n=3$ & $n=4$ \\ \hline
5 & (3,2) & - & - \\
7 & (5,2), (2,3) & - & - \\
11 & (9,2), (4,3) & (3,2,2), (2,3,2) & - \\
13 & (11,2), (5,3), (3,4), (2,5) & - & - \\
17 & (15,2), (7,3), (3,5) & (5,2,2), (2,5,2) & - \\
19 & (17,2), (8,3), (5,4), (2,7) & (4,3,2), (3,4,2), (3,2,3) & - \\
23 & (21,2), (10,3) & (7,2,2), (5,3,2), (3,5,2), (2,7,2) & - \\
29 & 4 tuples & 6 tuples & (3,2,3,2) \\
31 & 6 tuples & 2 tuples & (5,2,2,2) 
\end{tabular}
\end{examples}

\medskip

{\bf Proof of theorem \ref{t6.1}:}
Let $\kappa:N\to N$ be a map as in section \ref{s3}, so for any $j\in N$
the sets $J=\{j\}$ and $K=\{\kappa(j)\}$ satisfy (C2) for $R=\supp f$.

The proof proceeds in 4 steps: Step 1 extends some notations and formulas
from the proof of theorem \ref{t4.3}. Step 2 shows that $\kappa$ is of chain type.
Step 3 shows all remaining statements in (a) and (b). Step 4 proves part (c).

\medskip
{\bf Step 1.}
We consider the graph of $\kappa$.
The union of components $C$ with $|C|=1$ is called $C_{Fermat}$.
For a component $C$ let $C_1\subset C$ be the 
root of $C$ if $C$ is a tree, the cycle in $C$ if $C$ contains a cycle, and
$C_1=C$ if $|C|=1$.

For a component $C$ with $|C|\geq 2$ let $C_2\subset C-C_1$ be the set of leaves,
that is, the vertices without incoming arrows, and let $C_3\subset C-C_2$ be 
the set of branch points, that is, the vertices with $\geq 2$ incoming arrows.
The multiplicity $r(j)\in \N$ of a branch point $j\in C_3$ is the number of 
incoming arrows minus $1$. If $C$ is not a cycle then $C_3\neq\emptyset$, 
$C_3\cap C_1\neq\emptyset$ and $\sum_{c\in C_3}r(c)=|C_2|$.

The union of all leaves is called $C_{Leaf}$, the union of all branch points 
is called $C_{Branch}$.

For a component $C$ with $|C|\geq 2$ and for $j\in C$ let
$$
\wwh C(j)=(j,\kappa(j),...,\kappa^{\wwh l(j)}(j))
$$
be the longest tuple witout repetition: If $C$ is a tree then $\kappa^{\wwh l(j)}(j)$
is the root and $\kappa^{\wwh l(j)-1}(j)$ is not the root.
If $C$ contains a cycle, $\wwh C(j)$ hits the cycle and runs around it almost once,
so it hits the cycle in $\kappa^{\wwh l(j)+1}(j)$.
If $k\in\wwh C(j)$ let $C(j,k)$ be the tuple from $j$ to $k$,
$C(j,k)=(j,\kappa(j),...,k)$.

The definition of $C(j)$ in the proof of theorem \ref{t4.3} is slightly changed here:
For $j\in C-C_1$ (not only $j\in C_2$), 
let $C(j)=(j,\kappa(j),...,\kappa^{l(j)}(j))$ be the sub-tuple
of $\wwh C(j)$ which stops just before reaching $C_1$, so
$\kappa^{l(j)}(j)\notin C_1$, $\kappa^{l(j)+1}(j)\in C_1$. 
For $j\in C_1$ define $C(j):=\emptyset$.

For any $C(j,k)$ define with $\rho$ as in \eqref{3.2}
$$
\wwh \rho(C(j,k)):= \rho(j,\kappa(j),...,k)
$$
and define $\wwh\rho(\emptyset):=1$.

Now formula \eqref{3.1} for the weight $w_j$ of a vertex $j\in C_1$ on a cycle
can be rephrased as 
\begin{eqnarray}\label{6.2}
w_j=\frac{\wwh \rho(\wwh C(j)-\{j\})}{\prod_{k\in C_1}a_k-(-1)^{|C_1|}}.
\end{eqnarray}
And formula \eqref{4.11} generalizes to
\begin{eqnarray}\label{6.3}
\mu(C(j,k))= \prod_{i\in C(j,k)} \left(\frac{1}{w_i}-1\right)
= \frac{\wwh \rho(C(j,k)) + (-1)^{|C(j,k)|+1} w_{\kappa(k)}}
{1-w_{\kappa(k)}}.
\end{eqnarray}
For $j,k\in C-C_1$ and $k\in C(j)$ the tuple $\wwh C(j)$ contains the tuple
$\wwh C(\kappa(k))$, they hit the cycle or root $C_1$ at the same vertex
$l_1\in C_1$ and end at the same vertex $l_2\in C_1$, with $\kappa(l_2)=l_1$.
For such $j$ and $k$ one calculates with \eqref{6.2} and \eqref{6.3}
\begin{eqnarray}\nonumber
&&\mu(C(j,k)) = \frac{\mu(\wwh C(j))}{\mu(\wwh C(\kappa(k)))}\\
&=& \frac{\wwh \rho(\wwh C(j)) + (-1)^{|\wwh C(j)|+1}w_{l_1}}
{\wwh \rho(\wwh C(\kappa(k))) + (-1)^{|\wwh C(\kappa(k))|+1}w_{l_1}}\nonumber\\
&=& 
\frac{\left(\prod_{l\in C_1}a_l-(-1)^{|C_1|}\right)\wwh\rho(\wwh C(j)) 
+ (-1)^{|\wwh C(j)|+1}\wwh \rho(\wwh C(l_1)-\{l_1\})}
{\left(\prod_{l\in C_1}a_l-(-1)^{|C_1|}\right)\wwh\rho(\wwh C(\kappa(k))) 
+ (-1)^{|\wwh C(\kappa(k))|+1}\wwh \rho(\wwh C(l_1)-\{l_1\})}\nonumber\\
&=& 
\frac{\left(\prod_{l\in C_1}a_l\right)\wwh\rho(\wwh C(j)) 
+ (-1)^{|C_1|+1}\left(\prod_{l\in C_1}a_l\right)\wwh \rho(C(j))}
{\left(\prod_{l\in C_1}a_l\right)\wwh\rho(\wwh C(\kappa(k))) 
+ (-1)^{|C_1|+1}\left(\prod_{l\in C_1}a_l\right)
\wwh \rho(C(\kappa(k)))}\nonumber\\
&=& 
\frac{\wwh\rho(\wwh C(j)) 
+ (-1)^{|C_1|+1}\wwh \rho(C(j))}
{\wwh\rho(\wwh C(\kappa(k))) 
+ (-1)^{|C_1|+1}\wwh \rho(C(\kappa(k))}\label{6.4}
\end{eqnarray}

A component $C$ with $|C|\geq 2$ which is not a cycle is a tree or a cycle with attached 
trees. One can choose a map $\beta:C_2\to C_3$ from the leaves to the branch points
such that $k\in C_3$ is the image of $r(k)$ leaves and $\beta(j)\in\wwh C(j)$ for any leaf $j$.
Then $C-C_1$ is the disjoint union 
$\bigcup_{j\in C_2}(C(j)-C(\beta(j)))$, here the sets underlying the tuples
are meant. Therefore
\begin{eqnarray}\nonumber
\mu(C)&=& \prod_{j\in C_1}a_j\cdot \prod_{j\in C_2}\mu(C(j))\cdot
\prod_{j\in C_3}\mu(C(j))^{-r(j)}\\
&=& \prod_{j\in C_1}a_j\cdot 
\prod_{j\in C_2}\left(\wwh \rho(\wwh C(j))+(-1)^{|C_1|+1}\wwh \rho(C(j))\right)\nonumber\\
&& \cdot 
\prod_{j\in C_3}\left(\wwh \rho(\wwh C(j))+(-1)^{|C_1|+1}\wwh \rho(C(j))\right)^{-r(j)}
\label{6.5}
\end{eqnarray}

\medskip
{\bf Step 2.}
If $C_{Leaf}=\emptyset$ then the graph of $\kappa$ is a union of points and cycles,
and 
\begin{eqnarray*}
\mu=\prod_{j\in C_{Fermat}}(a_j+1)\cdot \prod_{C\ cycle}\prod_{j\in C}a_j.
\end{eqnarray*}
Then $\mu=$ prime number and all $a_j\geq 2$ imply $n=1$.

So suppose $C_{Leaf}\neq\emptyset$. Then there is a leaf $j_0\in C_{Leaf}$ such that
compared to all leaves $j\in C_{Leaf}$ the number
$\wwh\rho(\wwh C(j))+(-1)^{|C_1|+1}\wwh \rho(C(j))$ is maximal for $j=j_0$.
Here and later by a slight abuse of notation we denote for any $j\in C-C_{Fermat}$ 
the cycle or root in the component of $j$ by $C_1$.
Now choose a map $\beta:C_{Leaf}\to C_{Branch}$ as at the end of step 1 and 
with the additional property $\beta(j_0)\in C_1$, so $\wwh C(j_0)$
hits $C_1$ in $\beta(j_0)$. This is possible. Define the following natural numbers
\begin{eqnarray*}
A_0&:=& \wwh\rho(\wwh C(j_0))+(-1)^{|C_1|+1}\wwh \rho(C(j_0)),\\
A_j&:=& a_j+1\qquad \textup{for }j\in C_{Fermat},\\
A_j&:=& a_j \qquad \textup{for }j\in \bigcup_{cycles\ C}C,\\
B_{j_0}&:=& \mu(C_1)=\prod_{j\in C_1}a_j \qquad (C_1\textup{ for }j_0),\\
B_j&:=& \wwh\rho(\wwh C(j))+(-1)^{|C_1|+1}\wwh \rho(C(j))\qquad
\textup{for }j\in C_{Leaf}-\{j_0\},\\
D_j&:=& \wwh\rho(\wwh C(\beta(j)))+(-1)^{|C_1|+1}\wwh \rho(C(\beta(j))) \qquad
\textup{for }j\in C_{Leaf}.
\end{eqnarray*}
Then 
\begin{eqnarray}
\mu&=& A_0\cdot\prod_{j\in C_{Fermat}\cup(\textup{all cycles})} A_j
\cdot \prod_{j\in C_{Leaf}} \frac{B_j}{D_j},\label{6.6}
\end{eqnarray}
and
\begin{eqnarray*}
&& \textup{all } A_j\geq 2,
\nonumber\\
&&A_0\geq B_j\qquad \textup{for }j\in C_{Leaf}-\{j_0\},\\
\frac{B_j}{D_j}&=& \mu(C(j,\beta(j))-\{\beta(j)\})>1\qquad\textup{for }
j\in C_{Leaf}-\{j_0\}.
\end{eqnarray*}
For $j=j_0$ the map $\beta$ was chosen with $\beta(j_0)\in C_1$, so 
$C(\beta(j_0))=\emptyset$, $\wwh\rho(C(\beta_0))=1$, and 
\begin{eqnarray*}
\frac{B_{j_0}}{D_{j_0}} = 
\frac{\prod_{j\in C_1}a_j}
{\wwh\rho(\wwh C(\beta(j_0)))+(-1)^{|C_1|+1}\wwh \rho(C(\beta(j_0)))}
\left\{\begin{matrix}=1\textup{ if }|C_1|=1\\
>1\textup{ if }|C_1|>1\end{matrix}\right.
\end{eqnarray*}
And
\begin{eqnarray*}
A_0 &=& \left(B_{j_0}-(-1)^{|C_1|}\right)\wwh\rho(C(j_0))
+ (-1)^{|C(j_0)|+1}\wwh \rho(\wwh C(\beta(j_0)-\{\beta(j_0)\})),\\
1 &\leq& \wwh \rho(\wwh C(\beta(j_0))-\{\beta(j_0)\}) < B_{j_0},\\
1&\leq & \wwh \rho (C(j_0))\quad \textup{ and }\quad\\
3&\leq & \wwh \rho (C(j_0))\textup{ if } |C(j_0)|\geq 2,
\end{eqnarray*}
so always
\begin{eqnarray*}
A_0\geq B_{j_0}.
\end{eqnarray*}
Summarizing, we obtain
\begin{eqnarray}\label{6.7}
&&A_j<\mu\textup{ for }j\neq 0,\quad B_j\leq A_0,\quad A_0\leq \mu,\\
A_0=\mu &\iff& 
C_{Fermat}\cup(cycles)=\emptyset,\quad C_{Leaf}=\{j_0\},\quad |C_1|=1\nonumber\\
&\iff& \kappa\textup{ is of chain type with the chain }\wwh C(j_0).\label{6.8}
\end{eqnarray}
$\mu$ is a prime number by assumption. It must divide one of the factors
$A_j$ or $B_j$ in \eqref{6.6}. Because of \eqref{6.7} this forces $A_0=\mu$.
Because of \eqref{6.8} $\kappa$ is of chain type with the chain $\wwh C(j_0)$.

\medskip
{\bf Step 3.}
After renumbering of the vertices of its graph, $\kappa:N\to N$ is the map with $
\kappa(1)=1$, $\kappa(i)=i-1$ for $2\leq i\leq n$.
Then $f$ contains the monomials $x_1^{a_1+1}, x_2^{a_2}x_1,...,x_n^{a_n}x_{n-1}$.
The Milnor number is
\begin{eqnarray*}
\mu = A_0=\wwh \rho(\wwh C(n)) + \wwh \rho(C(n)) = \rho(a_n,...,a_2,a_1+1).
\end{eqnarray*}
The weights $w_i$ and the numbers $s_i,t_i\in\N$ with $gcd(s_i,t_i)=1$, $w_i=\frac{s_i}{t_i}$
are determined recursively by 
$w_1=\frac{1}{a_1+1}$, $s_1=1$, $t_1=a_1+1$,
\begin{eqnarray*}
\frac{s_{i+1}}{t_{i+1}}&=&w_{i+1}=\frac{1-w_i}{a_{i+1}} = \frac{t_i-s_i}{t_i\cdot a_{i+1}},\\
s_{i+1}&=&\frac{t_i-s_i}{\gamma_i},\ t_{i+1}=\beta_it_i,\\
\textup{where }&&a_{i+1}=\beta_i\gamma_i,\ \gamma_i=\gcd(a_{i+1},t_i-s_i).
\end{eqnarray*}
Thus
\begin{eqnarray*}
\mu = \prod_{i=1}^n \left(\frac{1}{w_i}-1\right) 
=\prod_{i=1}^n\frac{t_i-s_i}{s_i}
=\gamma_1\cdot ...\cdot \gamma_{n-1}
\cdot (t_n-s_n).
\end{eqnarray*}
$\mu$ being a prime number forces $\gamma_i=1$, $\beta_i=a_{i+1}$, $s_{i+1}=t_i-s_i$
and 
\begin{eqnarray*}
t_i&=&a_it_{i-1}=a_i...a_2\cdot (a_1+1),\\
s_i &=&\rho(a_{i-1},...,a_2,a_1+1)=t_{i-1}-t_{i-2}+...
+(-1)^{i-1}.
\end{eqnarray*}

Finally we show that the only monomials of weighted degree $d$ are 
$x_1^{a_1+1}, x_2^{a_2}x_1,...,x_n^{a_n}x_{n-1}$.
Then $f$ is as claimed in (a).
Let $\sum_{i=1}^n\delta_ie_i\in (\N_0^n)_d$. 
Let $j$ be maximal with $\delta_j>0$. Then
\begin{eqnarray*}
\delta_j\cdot\frac{s_j}{t_j} = 1-\sum_{i<j}\delta_i\cdot \frac{s_i}{t_i}.
\end{eqnarray*}
The denominator of the rational number on the right hand side is a divisor of $t_{j-1}$,
and $t_j=a_jt_{j-1}$. Therefore $\delta_j=a_j\varepsilon$ for some $\varepsilon\in\N$.
But
\begin{eqnarray*}
a_jw_j+w_{j-1}=1,\quad \textup{so }2a_jw_j>1,\quad\textup{so }\varepsilon=1,\quad
\textup{so }\sum_{i<j}\delta_iw_i=w_{j-1}.
\end{eqnarray*}
Then $\delta_{j-1}=1$, $\delta_i=0$ for $i<j-1$, so $\sum_i\delta_ie_i=a_je_j+e_{j-1}$.

\medskip
{\bf Step 4.}
Following \cite{MO}, we define the divisor $\divis p(t)$ of a unitary polynomial
$p(t)=\prod_{i=1}^k(t-\lambda_i)$ with zeros $\lambda_i\in S^1$ as the element
\begin{eqnarray*}
\divis p(t) := \sum_{i=1}^k \langle \lambda_j\rangle \in \Q[S^1]
\end{eqnarray*}
in the group ring $\Q(S^1)$. Denote $\Lambda_k:=\divis (t^k-1)$.
Then $1=\Lambda_1$ is a unit element and 
$\Lambda_a\cdot \Lambda_b = \gcd(a,b)\cdot \Lambda_{lcm(a,b)}$.

By \cite[Theorem 4]{MO} the divisor of the characteristic polynomial $\Delta(t)$ 
of the monodromy of $f$ is
\begin{eqnarray*}
\divis\Delta(t) =\prod_{i=1}^n\left(\frac{1}{s_i}\Lambda_{t_i}-1\right).
\end{eqnarray*}
Using $s_{i+1}=t_i-t_{i-1}+...+(-1)^i$ and 
$\Lambda_{t_i}\cdot \Lambda_{t_j}=t_i\cdot \Lambda_{t_j}$ for $i\leq j$,
we calculate
\begin{eqnarray*}
\divis\Delta(t) &=& (\Lambda_{t_1}-1)(\frac{1}{s_2}\Lambda_{t_2}-1)\cdot ... \\
&=& (\frac{t_1-1}{s_2}\Lambda_{t_2}-\Lambda_{t_1}+1)\cdot ...
= (\Lambda_2-\Lambda_1+1)\cdot ...\\
&=& ... = \Lambda_{t_n}-\Lambda_{t_{n-1}}+...+(-1)^{n-1}\Lambda_{t_1}+(-1)^n.
\end{eqnarray*}
This shows part (c) of theorem \ref{t6.1}
\hfill$\Box$

\bigskip
For a fixed $n\in \N$ a natural number $\mu>2^n$ is called an {\it $n$-gap}
if there does not exist a quasihomogeneous polynomial 
$f\in{\bf m}^3_{\C^n,0}$ with an isolated singularity at $0$ and
Milnor number $\mu$. 

\begin{corollary}\label{t6.3}
For $n\geq 3$ the set of $n$-gaps contains the set
\begin{eqnarray*}
\{2p+(-1)^n\, |\, p\textup{ and } 
2p+(-1)^n\textup{ are prime numbers},2p+(-1)^n>2^n\}.
\end{eqnarray*}
\end{corollary}

{\bf Proof:}
Consider a $p\in\N$ such that $\mu=2p+(-1)^n$ is bigger than $2^n$ 
and is a prime number, but not an $n$-gap. Then by theorem \ref{t6.1}
there exist $a_1,...,a_n\in\N-\{1\}$ with
\begin{eqnarray*}
2p+(-1)^n &=& \rho(a_n,...,a_2,a_1+1) \\
\textup{thus  }\quad 2p&=& (a_1+1)(\rho(a_n,...,a_2)+(-1)^{n-1}).
\end{eqnarray*}
But $a_1+1\geq 3$ and $\rho(a_n,...,a_2)+(-1)^{n-1}\geq 3$ if $n\geq 3$,
thus $p$ cannot be a prime number.
\hfill $\Box$

\begin{remarks}\label{t6.4}
(i) \cite{Ri} A natural number $p$ such that $p$ and $2p+1$ are prime numbers
is called a {\it Sophie Germain prime number}. There are conjectures
of Dickson (1904) (and a generalization called 
{\it hypothesis H} of Schinzel (1956)) and of Hardy and Littlewood (1923)
which would imply that the set of Sophie Germain prime numbers 
as well as the set $\{p\, |\, p\textup{ and }2p-1\textup{ are prime numbers}\}$
are infinite. But the infinity of both sets seems to be unknown.

(ii) It is also interesting to ask how many other $n$-gaps exist for $n\geq 3$.
There are 20 $3$-gaps with $8<\mu\leq 1000$, 19 of them are of the type $2p-1$
with $p$ and $2p-1$ being prime numbers, $9$ is the only other gap.
There are 21 $4$-gaps with $16<\mu\leq 500$, 14 of them are of the type
$2p+1$ with $p$ a Sophie Germain prime number, the other ones are 
17, 18, 19, 27, 74, 219, 314.
\end{remarks}

\end{document}